\documentclass[preprint,11pt]{elsarticle}

\usepackage{amssymb}

\usepackage{color}
\usepackage{amsmath}
\usepackage{comment}

\newcommand{\bx}{\mathbf{x}}
\newcommand{\bv}{\mathbf{v}}
\newcommand{\bB}{\mathbf{B}}

\newcommand{\bE}{\mathbf{E}}
\newcommand{\bF}{\mathbf{F}}

\newcommand{\bj}{\mathbf{j}}

\newcommand{\Dt}{\Delta t}

\newcommand{\nph}{^{n+1/2}}

\newcommand{\order}[1]{\mathcal{O} \left( #1 \right)}

\newcommand{\rd}{\,\mathrm{d}}

\DeclareMathOperator*{\argmin}{arg\,min}

\journal{Journal of Computational Physics}

\begin{document}

\begin{frontmatter}



\title{An explicit, energy-conserving particle-in-cell scheme}


\author[LLNL]{Lee F.\ Ricketson}

\affiliation[LLNL]{organization={Lawrence Livermore National Laboratory, Center for Applied Scientific Computing},
            addressline={7000 East Avenue}, 
            city={Livermore},
            postcode={94550}, 
            state={CA},
            country={USA}}

\author[UW]{Jingwei Hu}

\affiliation[UW]{organization={Department of Applied Mathematics, University of Washington},
            addressline={Box 353925}, 
            city={Seattle},
            postcode={98195}, 
            state={WA},
            country={USA}}

\begin{abstract}
We present an explicit temporal discretization of particle-in-cell schemes for the Vlasov equation that results in exact energy conservation when combined with an appropriate spatial discretization.  The scheme is inspired by a simple, second-order explicit scheme that conserves energy exactly in the Eulerian context.  We show that direct translation to particle-in-cell does not result in strict conservation, but derive a simple correction based on an analytically solvable optimization problem that recovers conservation.  While this optimization problem is not guaranteed to have a real solution for every particle, we provide a correction that makes imaginary values extremely rare and still admits $\mathcal{O}(10^{-12})$ fractional errors in energy for practical simulation parameters.  We present the scheme in both electrostatic -- where we use the Amp\`{e}re formulation -- and electromagnetic contexts.  With an electromagnetic field solve, the field update is most naturally linearly implicit, but the more computationally intensive particle update remains fully explicit.  We also show how the scheme can be extended to use the fully explicit leapfrog and pseudospectral analytic time-domain (PSATD) field solvers.  The scheme is tested on standard kinetic plasma problems, confirming its conservation properties.   
\end{abstract}



\begin{keyword}
particle-in-cell \sep energy conservation \sep plasma \sep Vlasov
\end{keyword}

\end{frontmatter}


\section{Introduction}
\label{sec:intro}
Since its inception over sixty years ago, the particle-in-cell (PIC) scheme has been an enormously popular method for solving the Vlasov equation and its related models of kinetic plasma dynamics.  PIC's desirable features include mitigation of the curse of dimensionality -- due to its lack of a mesh in velocity space -- simplicity, robustness, and scalability.  PIC codes are routinely used to simulate complex plasma physical phenomena on large-scale supercomputers \cite{ku2009full, fiuza2011three, chen2007electromagnetic, vay2021modeling}. 

However, PIC is not without drawbacks.  In addition to the obvious consideration of particle sampling noise that can pollute solutions, the finite grid instability \cite{birdsall2018plasma} has long plagued PIC schemes.  The most commonly observed symptom of the finite grid instability is so-called ``grid heating", which causes plasmas simulated by PIC feature a secular growth in total energy over time.  As a result, there has been considerable interest in \textit{energy conserving} PIC schemes for many years.  

Much of this effort has focused on \textit{implicit} schemes.  This began in the 1980s, when predictor-corrector style schemes that mimicked implicitness were developed \cite{friedman1981direct, langdon1983direct}.  These schemes improved energy conservation relative to explicit schemes, but the conservation was not exact and they featured limited improvements to stability.  In the 2010s, \textit{fully implicit} PIC schemes were developed by Chen, Chac\'{o}n and collaborators \cite{chen2011energy, chen2015multi, chacon2016curvilinear}.  These schemes were enabled by modern linear and nonlinear solvers, working in concert with a clever decomposition of the nonlinear system resulting from the implicit discretization.  Among many other desirable features, these schemes feature \textit{exact} energy conservation (in practice, conservation is achieved up to solver tolerances), and have been shown to be more resistant to the finite grid instability than explicit PIC schemes \cite{barnes2021finite}.  Around the same time, Lapenta and collaborators developed \textit{semi-implicit} PIC schemes that also conserve energy exactly \cite{markidis2011energy, lapenta2017exactly, bacchini2019relativistic}.  

While the fully- and semi-implicit schemes have several distinguishing characteristics, we note here only that the semi-implicit schemes treat the particle update explicitly while fully implicit schemes treat the particle push implicitly.  In general, this means that semi-implicit schemes are cheaper on a per-timestep basis, but fully implicit schemes remain accurate and stable at larger time-steps.  Also noteworthy and relevant to our development is a different semi-implicit PIC scheme of Chen and Chac\'{o}n \cite{chen2020semi}.  There, in contrast to the semi-implicit work of Lapenta, the particle update is treated implicitly while the field solve is explicit.  

As hinted at above, implicit schemes come with the added benefit of allowing PIC to step over stiff time-scales in the problem.  The mitigation of the finite grid instability enjoyed by such schemes also permits the use of spatial cells much larger than the Debye length -- a tremendous advantage when solution structures can be adequately resolved on such meshes.  However, implicitness comes with a cost as well.  The solution of linear and -- for fully implicit schemes -- nonlinear systems can be computationally intensive, and the efficient, robust implementation of such solvers can be quite complex.  

As such, for plasmas systems featuring Debye-scale spatial structures and without stiff time-scales one wishes to step over, there is motivation to pursue an \textit{explicit} scheme that nevertheless features exact energy conservation.  Such a scheme would avoid the overhead of implicit solvers while still eliminating the long-standing problem of grid heating.  Such a scheme is the topic of this article.  

There have been two other recent efforts in this direction.  In \cite{gonoskov2024explicit}, a novel splitting scheme for the relativistic Vlasov-Maxwell system was introduced that results in an explicit scheme with exact energy conservation.  However, the nature of the splitting requires that many particle-based operations that could otherwise be parallelized be performed in serial.  This has severe negative consequences for the scalability of the algorithm.  More similar to our approach is that of \cite{ji2023asymptotic}, in which a global Lagrange multiplier approach is applied at the end of each time-step to recover energy conservation.  In addition to the scalability consequences of such a global operation, the impact this operation has on the order of accuracy of the scheme has not been studied.

Our approach applies a Lagrange multiplier optimization procedure to each particle \textit{individually}.  The locality of the approach means both that it has essentially no impact on parallel scalability and that we are able to analyze the accuracy of the scheme at the level of the particle update.  We show that the optimization problem that enforces energy conservation has an analytic solution and does not impact the second-order temporal accuracy of the scheme.  We also show that, like the commonly used Boris integrator, the magnetic field does identically zero work, which improves long-time accuracy over schemes without this feature.  Finally, we show that the scheme is compatible with exact charge conservation.  

The Lagrange multiplier optimization procedure does admit imaginary solutions, but we show that such solutions are exceedingly rare and propose a correction scheme that makes them even more rare so that energy errors approach double-precision round-off errors in our tests.  We also show how to avoid imaginary velocities in these rare cases.  

We derive our scheme in both the electrostatic and electromagnetic contexts.  In the electromagnetic case, it is well-known that in addition to energy conservation, accurate representation of the light-wave dispersion relation is important for long-time accuracy.  We thus show that our new temporal discretization is compatible with two spatial discretizations of Maxwell's equations that are commonly used and known to have good dispersion properties: the Yee lattice \cite{Yee66} and the pseudo-spectral analytic time domain (PSATD) method \cite{vay2013domain, lehe2018review}. The resulting scheme has the same essential structure as most existing explicit PIC methods, which should make it readily implementable in existing explicit PIC codes with only minor modifications.  

The remainder of the article is structured as follows.  In Section \ref{sec:bg}, we review the Vlasov equation, Maxwell's equations, PIC discretizations of each and the important properties of such discretizations.  In Section \ref{sec:main} we derive our method in the electrostatic case, show that it is second-order accurate in time, derive a correction procedure that reduces the probability of imaginary solutions to the optimization problem, describe the most straightforward extension of the electromagnetics, and show how charge conservation can be achieved.  In Section 4, we show that the scheme is compatible and retains energy conservation with the Yee lattice combined with a leapfrog-type temporal discretization of Maxwell's equation (together often called the FDTD method) as well as the PSATD Maxwell discretization.  In Section \ref{sec:numerics}, we apply the scheme to linear Landau damping, the two-stream instability and the Weibel instability to confirm the scheme's conservation properties.  We conclude and describe directions for future work in Section \ref{sec:conc}.

\section{Background}
\label{sec:bg}

\subsection{Vlasov Equation}
We are concerned with the numerical solution of the Vlasov equation:
\begin{equation} \label{eq:vlasov}
 \frac{\partial f}{\partial t}  + \bv \cdot \nabla_{\bx} f + \left( \bE + \bv \times \bB \right) \cdot \nabla_{\bv} f = 0,
\end{equation}
where $f(\bx, \bv, t)$ typically represents the phase-space number density of some species of charged particle.  Here and throughout the remainder of the paper, we work in a dimensionless formulation in which time has been scaled by the plasma frequency $\omega_p = \sqrt{n_0 e^2 / m \epsilon_0}$, space by Debye length $\lambda_D = \sqrt{\epsilon_0 T / n_0 e^2}$,
 and velocity by a thermal speed $v_{th} = \sqrt{T/m} = \omega_p \lambda_D$.  Here, $n_0$ is some reference number density, $m$ the species mass, $e$ the fundamental unit charge, $\epsilon_0$ the permittivity of free space, and $T$ a reference temperature.

The equation is closed by expressing the electromagnetic fields $\bE(\bx, t)$ and $\bB(\bx, t)$ in terms of $f$.  Two common forms of this closure are Maxwell's equations and their electrostatic approximation.  In our dimensionless formulation, Maxwell's equations read
\begin{align} \label{eq:maxwell}
    \nabla_{\bx} \cdot \bE &= \rho - 1, \\
    \nabla_{\bx} \cdot \bB &= 0, \\
    \nabla_{\bx} \times \bE &= -\frac{\partial \bB}{\partial t}, \\
    \nabla_{\bx} \times \bB &= \frac{1}{c^2} \left( \frac{\partial \bE}{\partial t} + \bj \right),
\end{align}
with charge density $\rho$ and current density $\bj$ given by
\begin{equation}
    \rho = \int f \rd\bv, \qquad \bj = \int \bv f \rd\bv,
\end{equation}
and the speed of light is again in our dimensionless variables -- that is, $c$ here denotes the physical speed of light divided by $v_{th}$.  We also note that we have included a uniform background charge density in Gauss' law, intended to capture an ion density that is constant on time-scales of interest.  Generalizations to multi-species systems are straightforward.  

The electrostatic approximation corresponds to the assumption that $\bE = -\nabla_{\bx} \phi$, and has two equivalent formulations in the continuum: the Poisson form
\begin{equation} \label{eq:poisson}
    -\Delta_{\bx} \phi = \rho - 1,
\end{equation}
and the Amp\`{e}re form
\begin{equation} \label{eq:ampere}
    \frac{\partial \Delta_{\bx} \phi}{\partial t}  = \nabla_{\bx} \cdot \mathbf{j}.
\end{equation}
Equivalence of these two forms is enforced by the continuity equation, derived by integrating the Vlasov equation \eqref{eq:vlasov} over all velocity space:
\begin{equation} \label{eq:continuity}
   \frac{\partial \rho}{\partial t}  + \nabla_{\bx} \cdot \bj = 0.
\end{equation}
Elementary substitution shows that, given \eqref{eq:continuity}, if \eqref{eq:poisson} is satisfied initially then \eqref{eq:ampere} ensures it is satisfied at all later times.  Clearly, the combination of \eqref{eq:poisson} and \eqref{eq:continuity} also implies \eqref{eq:ampere}.  The continuity equation plays an analogous role in the electromagnetic (i.e., Maxwell) context: it ensures that the first two lines of \eqref{eq:maxwell}-(5) are satisfied for all time if they are satisfied initially.  These equations are thus often called ``involutions", and not explicitly enforced (except at $t=0$) in many discretizations.  

\subsection{Particle-in-Cell Discretization}
In the most general case, $\bx, \bv \in \mathbb{R}^3$, making the curse of dimensionality weigh quite heavily on numerical simulation of the Vlasov equation.  As a result, particle-based methods have long been popular for its solution, by far the most common of which being the particle-in-cell (PIC) schemes.  These begin with the ansatz that $f$ may be approximated by a sum of $N_p$ Dirac delta functions in phase space: 
\begin{equation}
    f = \sum_{p=1}^{N_p} w_p \delta(\bx - \bx_p(t)) \delta(\bv - \bv_p(t)).
\end{equation}
The evolution equations for the locations and magnitudes of these delta functions are clearly the characteristic equations:
\begin{equation} \label{eq:characteristicODE}
    \frac{\rd\bx_p}{\rd t} = \bv_p, \qquad \frac{\rd \bv_p}{\rd t} = \bE(\bx_p, t) + \bv_p \times \bB(\bx_p, t), \qquad \frac{\rd w_p}{\rd t} = 0.
\end{equation}
The system is closed by introducing a mesh in configuration space, whose grid points we index by $h$, and approximating $\rho$ and $\bj$ on that mesh:
\begin{equation}
    \rho_h(t) = \frac{1}{| \mathbf{h} |} \sum_p w_p S^h_\rho (\bx_h - \bx_p(t)), \qquad \bj_h(t) = \frac{1}{| \mathbf{h} |} \sum_p w_p \bv_p(t) S^h_\bj(\bx_h - \bx_p(t)),
\end{equation}
where $\bx_h$ is a point on the mesh and $| \mathbf{h} |$ denotes cell volume. $S^h_\rho/|\mathbf{h}|$ and $S^h_\bj/|\mathbf{h}|$ are approximate delta functions traditionally called ``shape functions" in the literature, owing to their alternative interpretation as the charge density associated with a single computational particle.  It is common to use $B$-splines as shape functions, with the most common being the first-order ``tent" function $S^h(z) = \max \{ 0, 1 - |z|/h \}$ (extension to multiple dimensions is achieved via tensor products).

The quantities $\rho_h$ and $\bj_h$ are used to solve for $\bE_h(t)$ and $\bB_h(t)$ on the mesh using Maxwell's equations or either of the equivalent electrostatic approximations described above and a spatial discretization of the user's choice.  Fields at the particle locations are recovered via interpolation, which may be written 
\begin{equation}
\begin{split}
    \bE(\bx_p,t) &= \sum_h \bE_h(t) S^h_\bE(\bx_p(t) - \bx_h), \\
    \bB(\bx_p,t) &= \sum_h \bB_h(t) S^h_\bB(\bx_p(t) - \bx_h).
\end{split}
\end{equation}
It is often the case that several of the shape functions $S^h_\rho$, $S^h_\bj$, $S^h_\bE$, $S^h_\bB$ are identical.  This may even be desirable in some contexts, but is not necessary, so we leave our notation general here.  A temporal discretization of the resulting system of ordinary differential equations completes the scheme.  

It is instructive to review some of the more commonly used temporal and spatial discretization choices because understanding what properties make them desirable informs the development of the scheme we present below.  We do so in the proceeding subsections, highlighting the desirable features of each that any new PIC scheme should seek to replicate.  


\subsubsection{Boris particle push} \label{sec:Boris}
By far the most common discretization of the characteristic equation \eqref{eq:characteristicODE} -- often called the ``particle push" step in the PIC literature -- is the so-called ``leapfrog" or ``Boris" method.  The method is usually presented with the $\bx$ and $\bv$ variables offset by a half step in time, but the analogy to our later development is closer in its equivalent ``symmetric" form, so that is the form we choose to write first:
\begin{equation} \label{eq:symm_Boris}
\begin{split}
    \bx_p^{n,*} &= \bx_p^n + \frac{\Delta t}{2} \bv_p^n, \\
    \bv_p^{n+1} &= \bv_p^n + \Delta t \left( \bE \left(\bx_p^{n,*}, t^{n+1/2} \right) + \bv_p^{n+1/2} \times \bB \left( \bx_p^{n,*}, t^{n+1/2} \right) \right), \\
    \bx_p^{n+1} &= \bx_p^{n,*} + \frac{\Delta t}{2} \bv_p^{n+1}.
\end{split}
\end{equation}
Here, $n$ indexes time-step and $\bv_p^{n+1/2} = (\bv_p^n + \bv_p^{n+1})/2$.  The more traditional time-offset form can be recovered by noting that $\bx_p^{n,*} = \bx_p^{n-1,*} + \Delta t \bv_p^n$ and defining $\bx_p^{n,*}$  
as the fundamental configuration-space unknown, often denoted as $\bx_p^{n+1/2}$.  However, in this paper we will reserve the $n+1/2$ superscript for quantities that are exact averages of quantities at steps $n$ and $n+1$ to avoid confusion in the later development.  Note also that while the scheme is technically implicit in the velocity variable, the implicitness is linear, the linear system to be solved has only three variables and is thus easily solved by hand.  The scheme is thus effectively explicit when it comes to computational efficiency.  

The Boris scheme preserves phase space volume \cite{qin2013boris} and has near conservation of energy over very long times in the special cases where the magnetic field is constant or the electric potential is quadratic \cite{HL18}.
This is at least partially because it preserves exactly in the discrete the fact that the magnetic field can do no work. 
The derivation is trivial: when $\bE = 0$, dot the velocity update with $\bv_p^{n+1/2}$ to find that $\| \bv_p^{n+1} \|^2 = \| \bv_p^n \|^2$.  This simple fact goes a long way toward explaining why the Boris scheme is preferred for long-time simulatons over higher-order numerical methods like RK4, as the higher order methods typically allow the magnetic field to do small but non-zero work which may accumulate into significant errors over time.  Keeping this property is highly desirable in any new temporal discretization.  


\subsubsection{Spatial discretizations and light-wave dispersion} 
\label{sec:space} For the Poisson form of the electrostatic field solve, all the usual methods of solving the Poisson equation may be used -- finite differences, finite elements, and pseudospectral methods have received particular attention.  The same may be said of the Amp\`{e}re form when combined with a standard temporal discretization.  The electromagnetic case, on the other hand, presents additional subtleties.  Particularly for problems in laser-plasma interaction, it is known to be important to choose a field discretization that accurately approximates the analytic dispersion relation for light waves in order to mitigate numerical Cherenkov radiation.  

The standard finite difference scheme uses the Yee's lattice \cite{Yee66}, in which components of the electric field are computed at the center of the cell edges to which they are parallel and magnetic field components on cell faces to which they are orthogonal.  This is combined with a leapfrog-style scheme in time, in which $\bE$ and $\bB$ are offset by half a temporal step.  Light wave dispersion is exact in one dimension with this method, with errors in multiple dimensions that are tolerable for many applications. 

More recently, pseudospectral methods for the Maxwell component of electromagnetic PIC schemes have received considerable attention \cite{lehe2018review}.  Of particular note is the so-called pseudospectral analytic time domain (PSATD) method \cite{vay2013domain}.  The spatial Fourier transform of the (non-involution) Maxwell's equations is 
\begin{equation}
\begin{split}
    \frac{\partial \mathcal{F}[\bB]}{\partial t} &= -i \mathbf{k} \times \mathcal{F}[\bE], \\
    \frac{\partial \mathcal{F}[\bE]}{\partial t} &= c^2 i \mathbf{k} \times \mathcal{F}[\bB] - \mathcal{F}[\bj],
\end{split}
\end{equation}
where $\mathcal{F}[\cdot]$ denotes the spatial Fourier transform.  If one fixes the current $\mathbf{j}$ at a specific time -- denote it $\tilde{t}$ -- this linear system of ODEs can be solved analytically over a single time-step.  That solution is (see Appendix A in \cite{vay2013domain})
\begin{equation} \label{eq:PSATDformula}
\begin{split}
    \mathcal{F} \left[ \bE^{n+1} \right] &= C \mathcal{F} \left[ \bE^n \right] + i S c \widehat{\mathbf{k}} \times \mathcal{F}\left[ \bB^n \right] - \frac{S}{kc} \mathcal{F} \left[ \bj \left( \tilde{t} \right) \right]\\
    & \qquad  + (1 - C) \widehat{\mathbf{k}} \left( \widehat{\mathbf{k}} \cdot \mathcal{F}\left[ \bE^n \right] \right) + \widehat{\mathbf{k}} \left( \widehat{\mathbf{k}} \cdot \mathcal{F}\left[ \bj(\tilde{t}) \right] \right) \left( \frac{S}{kc} - \Delta t \right), \\
    \mathcal{F}\left[ \bB^{n+1} \right] &= C \mathcal{F}\left[ \bB^{n} \right] - i\frac{S}{c} \widehat{\mathbf{k}} \times \mathcal{F}\left[ \bE^{n} \right] + i \frac{1-C}{kc^2} \widehat{\mathbf{k}} \times \mathcal{F} \left[ \bj\left( \tilde{t} \right) \right],
\end{split}
\end{equation}
where $C = \cos (kc\Delta t)$ and $S = \sin (kc\Delta t)$.  Typically one chooses $\tilde{t} = t^{n+1/2}$ or some approximation thereof to achieve second-order accuracy.  Note that the expression here differs from that in \cite{vay2013domain} only because of the different normalizations of Maxwell's equations -- in particular, our magnetic field $\bB$ differs from the one there by a factor of $c$.  

When preceded by a fast Fourier transform (FFT) and proceeded by an inverse FFT, this leads to a field solve with the exact correct light wave dispersion relation, as all vacuum terms are treated analytically and approximation is only introduced in the treatment of the current $\bj$.  

For our later analysis, it will be critical to realize that this formula more generally gives an expression for $\bE$ and $\bB$ at an arbitrary time $t^n+ \tau$, $\tau \in [0, \Dt]$.  This is achieved simply by letting $\Delta t \rightarrow \tau$ in the formula above.  The resulting expressions for $\bE$ and $\bB$ give analytic solutions to the ODE system
\begin{equation} \label{eq:PSATDconstruction}
\begin{split}
    \frac{\partial \bE_h}{\partial t} &= c^2 \nabla_h \times \bB_h - \mathbf{j}_h^{n+1/2}, \\
    \frac{\partial \bB_h}{\partial t} &= - \nabla_h \times \bE_h,
\end{split}
\end{equation}
with initial conditions given at time $t^n$ and the discrete curl operator defined in the pseudospectral manner.  


\subsubsection{Implicit PIC and Energy Conservation}
We also must note the development of implicit PIC schemes over the last decade.  Such methods feature exact (up to the tolerance of the linear and/or nonlinear solvers being employed) energy conservation, in contrast to their explicit counterparts.  They have also been shown to be more resistant to the finite grid instability than standard explicit schemes.  There are many variants of such schemes, and we discuss only a small subset here.  Following \cite{chen2011energy}, a fully implicit scheme in the electrostatic case has the form 
\begin{equation} \label{eq:implicitPIC}
\begin{split}
	\bx^{n+1}_p &= \bx^n_p + \Delta t\bv_p\nph, \\
	\bv^{n+1}_p &= \bv^n_p + \Delta t \left[ \bE\nph_p + \bv\nph_p \times \bB\nph_p \right], \\
	\bE_p\nph &= \sum_h \bE\nph_h S^h \left( \bx_p\nph - \bx_h \right), \\
	\bE_h^{n} &= -\nabla_h \phi_h^n, \\
	\nabla_h^2 \phi_h^{n+1} &= \nabla_h^2 \phi_h^n + \Dt \nabla_h \cdot \bj^{n+1/2}_h, \\
	\bj_h^{n+1/2} &= \frac{1}{ \lvert \mathbf{h} \rvert } \sum_p w_p \bv_p\nph S^h \left( \bx_h - \bx_p\nph \right).
\end{split}
\end{equation}
In the above, an $h$ subscript on a differential operator denotes its discretization on the mesh.  For generality, we do not specify this discretization here.  In keeping with the convention described above, all quantities at half-steps are averages of the corresponding values at neighboring integer steps.  $\bB$ is assumed to be some given analytic function, and $\bB_p^{n+1/2} = \bB(\bx_p^{n+1/2}, t^{n+1/2})$.  Note also that the shape functions used to interpolate the electric field and to deposit current are identical.  As a final note, observe that the implicit particle push (first two lines of \eqref{eq:implicitPIC}) shares the desirable feature of Boris that the magnetic field does identically zero work.  

The proof of energy conservation for such a scheme is by quite standard at this point.  Importantly for the development here, the only constraint placed on the spatial discretization is that two integration-by-parts identities must be respected in the discrete:
\begin{equation} \label{eq:IBP_ES}
\begin{split}
    |\mathbf{h}| \sum_h \bF_h \cdot \nabla_h f_h &= -|\mathbf{h}| \sum_h f_h \nabla_h \cdot \bF_h, \\
    |\mathbf{h}| \sum_h \nabla_h f_h \cdot \nabla_h g_h &= -|\mathbf{h}| \sum_h f_h \nabla^2_h g_h,
\end{split}
\end{equation}
for arbitrary grid quantities $f$, $g$ and $\bF$.  While this certainly doesn't hold for arbitrary spatial discretizations, it has been shown to hold for second-order finite differences on Yee's lattice \cite{chen2011energy} and pseudospectral methods \cite{ricketson2023pseudospectral} (although binomial filtering is required in the latter case under certain circumstances).  When not specifying a particular spatial discretization scheme, we will assume whatever scheme is chosen satisfies this identity.

This proof has been extended to the electromagnetic context.  Of particular note is the fact that a semi-implicit scheme can be derived by combining the implicit particle push above with Yee's leapfrog algorithm for the field update \cite{chen2020semi}, which is explicit and offsets $\bE$ and $\bB$ by half a time-step.  Exact energy conservation is achieved with the following non-standard definition of magnetic potential energy:
\begin{equation}
    W_B^n = \frac{1}{2} | \mathbf{h}| \sum_h  \bB_h^{n-1/2} \cdot \bB_h^{n+1/2}. 
\end{equation}
This definition is a second-order approximation of the more natural definition: $| \mathbf{h} | \sum_h \| \bB_h^n \|^2 / 2$.  Additionally, it is shown in the appendix of \cite{chen2020semi} that the non-standard definition is almost surely well-posed (i.e. non-negative).  


\subsubsection{Charge conservation}
\label{sec:chrgcons_bckgrnd}
Aside from the evolution-type Maxwell's equations, some care is required to ensure the involution equations are satisfied for all time at the discrete level.  Ensuring $\bB$ remains divergence-free is relatively straightforward: one just requires $\nabla_h \cdot (\nabla_h \times \mathbf{F}_h) = 0$ for arbitrary grid functions $\mathbf{F}_h$.  This is satisfied by numerous spatial discretizations, including the Yee's lattice and pseudospectral methods described above.  

Less trivial is enforcement of Gauss' law, $\nabla_h \cdot \bE = \rho - 1$.  Recall that this is ensured in the continuum by the combination of Amp\`{e}re's law and the continuity equation \eqref{eq:continuity}.  Because the charge density $\rho$ plays no role in the time-advancement of the system, one can simply \textit{define} $\rho$ in terms of some discretized version of the continuity equation \cite{ricketson2023pseudospectral, sturdevant2022eliminating}.  

While such a $\rho$ certainly guarantees Gauss' law is satisfied for all time (given a satisfactory initial $\rho$), there can be concern that such a $\rho$ can deviate over time from any reasonable definition based directly on the particles.  It is often demanded that there exists a well-defined shape function $S^h_\rho$ such that the definitions
\begin{equation}
    \rho_h^n = \frac{1}{|\mathbf{h}|} \sum_p w_p S^h_\rho(\bx_h - \bx_p^n), \qquad D_t \rho_h^n + \nabla_h \cdot \bj_h^{n+\alpha} = 0
\end{equation}
are \textit{equivalent} given appropriate initial data.  Here, $D_t$ is some discretization of the temporal derivative, and the parameter $\alpha$ indicates that current need not be evaluated at an integer time-step.  We say a shape function $S^h_\rho$ is well-defined if it is symmetric and $S^h_\rho(z)/h = \delta(z) + \mathcal{O}(h^2)$.  Failure to enforce such an equivalence has been shown to have deleterious effects in some circumstances \cite{chen2011energy}.  We shall call the definitions of $\rho$ above the \textit{direct deposition} and \textit{continuity} definitions, respectively.

The equivalence described above has been established in a variety of contexts.  The rough steps in most derivations may be described as follows.  By substituting the definition of current and the direct deposition definition of $\rho$ into the continuity definition, one finds that a sufficient condition for this equivalence is
\begin{equation}
    D_t S^h_\rho (\bx_h - \bx_p^n) + \nabla_h \cdot \left[ \bv_p^{n+\alpha} S^h_\bj (\bx_h - \bx_p^{n+\alpha}) \right] = 0
\end{equation}
for every particle $p$.  Guaranteeing this identity requires some smoothness of the shape functions $S^h_\bj$ and $S^h_\rho$ that is achievable within a given cell but not across cell boundaries.  One thus decomposes particle trajectories within a time-step into sub-steps, each of which lies within a single cell \cite{chen2011energy, chen2020semi, esirkepov2001exact}.  Current deposition and field interpolation are then performed on these sub-steps.  Additional details can be found in the references above, and we also describe the specific procedure used to achieve charge conservation in our scheme -- which follows these general steps -- in Section \ref{sec:chargeconservation}.


\section{The Method}
\label{sec:main}

\subsection{Motivating Eulerian Method}
To motivate our particle method, we first present a second-order energy-conserving Eulerian time discretization:
\begin{subequations} \label{eq:eulerian_scheme}
\begin{align}
&\frac{f^*-f^n}{\Delta t/2}+\bv \cdot \nabla_{\bx}f^n + (\bE^n+\bv \times \bB^n) \cdot \nabla_{\bv} f^n=0,\\
&\frac{\bE^{n+1}-\bE^n}{\Delta t}=c^2\nabla_\bx\times \bB^{n+1/2}-\bj^*, \label{Eulerian2}\\
&\frac{\bB^{n+1}-\bB^n}{\Delta t}=-\nabla_\bx\times \bE^{n+1/2},\label{Eulerian3}\\
&\frac{f^{n+1}-f^n}{\Delta t}+\bv \cdot \nabla_{\bx}f^* + (\bE^{n+1/2}+\bv \times \bB^{n+1/2}) \cdot \nabla_\bv f^*=0,\label{Eulerian4}
\end{align}
\end{subequations}
where 
$$
\bB^{n+1/2}=\frac{\bB^n+\bB^{n+1}}{2}, \quad
\bE^{n+1/2}=\frac{\bE^n+\bE^{n+1}}{2}, \quad
\bj^*=\int  \bv f^*\rd{\bv}.
$$
First note that the implicitness only appears in the field equations. In fact, the scheme is fully explicit if the magnetic field is zero. 

To see the energy conservation, multiplying \eqref{Eulerian4} by $\frac{1}{2}|\bv|^2$ and integrating in $\bv$ yields
\begin{equation} \label{moment1}
\frac{1}{\Delta t}\left(\int \frac{1}{2}|\bv|^2 f^{n+1} \rd{\bv}-\int \frac{1}{2}|\bv|^2 f^n \rd{\bv}\right)+\nabla_{\bx}\cdot\int \bv \frac{1}{2}|\bv|^2f\rd{\bv}= \bj^{*}\cdot\bE^{n+1/2}.
\end{equation}
On the other hand, multiplying \eqref{Eulerian2} by $\bE^{n+1/2}$ and \eqref{Eulerian3} by $\bB^{n+1/2}$ and adding them together, we obtain
\begin{equation} \label{moment2}
\frac{1}{2}\frac{|\bE^{n+1}|^2-|\bE^{n}|^2}{\Delta t}+\frac{c^2}{2}\frac{|\bB^{n+1}|^2-|\bB^{n}|^2}{\Delta t}
=c^2 \nabla_{\bx}\cdot\left(\bB^{n+1/2}\times \bE^{n+1/2}\right)-\bj^{*}\cdot\bE^{n+1/2}.
\end{equation}
Finally, integrating both \eqref{moment1} and \eqref{moment2} in $\bx$ (assuming periodic boundary condition in $\bx$) and adding them together gives the energy conservation:
\begin{align}
&\iint \frac{1}{2}|\bv|^2 f^{n+1} \rd{\bv}\rd{\bx}+\int \left(\frac{1}{2}|\bE^{n+1}|^2+\frac{c^2}{2}|\bB^{n+1}|^2\right)\rd{\bx}\nonumber\\
&=\iint \frac{1}{2}|\bv|^2 f^{n} \rd{\bv}\rd{\bx}+\int \left(\frac{1}{2}|\bE^{n}|^2+\frac{c^2}{2}|\bB^{n}|^2\right)\rd{\bx}.
\end{align}
This scheme is quite simple, and we are thus motivated to ask whether its translation to PIC can lead to a similarly simple particle scheme with exact energy conservation.  We show below that the answer is yes, but that an additional trick is required in the particle context.  

\subsection{PIC translation in electrostatic case}
\label{sec:initialscheme}
We begin in the simplest context: the Amp\`{e}re formulation of the electrostatic approximation with externally imposed magnetic field.  A relatively natural translation of the Eulerian scheme to the PIC context is as follows:
\begin{equation} \label{eq:naivetranslation}
\begin{split}
    \bx^*_p &= \bx^n_p + \frac{\Delta t}{2} \bv^n_p, \\
    \bv^*_p &= \bv^n_p + \frac{\Delta t}{2} \left( \bE^{n,*}_p + \bv_p^* \times \bB(\bx_p^*, t^{n+1/2}) \right), \\
    \bx^{n+1}_p &= \bx^n_p + \Delta t \bv^*_p, \\
	\nabla_h^2 \phi_h^{n+1} &= \nabla_h^2 \phi_h^n + \Dt \nabla_h \cdot \bj^{*,n+1/2}_h, \\
    \bv^{n+1}_p &= \bv^n_p + \Delta t \left( \bE^{n+1/2}_p + \bv_p^* \times \bB(\bx_p^*, t^{n+1/2}) \right), \\
\end{split}
\end{equation}
where we define
\begin{equation}
\begin{split}
    \bE^{n,*}_p &= \sum_h \bE_h^n S^h \left(\bx_p^* - \bx_h \right), \\
    \bj_h^{*,n+1/2} &= \frac{1}{| \mathbf{h} |} \sum_p w_p \bv_p^* S^h \left( \bx_p\nph - \bx_h \right),\\
     \bE^{n+1/2}_p &= \sum_h \bE^{n+1/2}_h S^h \left( \bx_p^{n+1/2} - \bx_h \right),
\end{split}
\end{equation}
and at any time level $\bE_h = -\nabla_h \phi_h$.  A slightly counter-intuitive choice made here is the definition of current, which uses $\bv_p^*$ as its velocity but $\bx_p^{n+1/2} = (\bx_p^{n+1} + \bx_p^n)/2$ as its location.  This choice will turn out to be critical for charge conservation later in the derivation.  The analysis can be simplified considerably by evaluating current and the electric field at $\bx_p^*$ instead of $\bx_p\nph$, but at the cost of charge conservation.  

Note that this scheme is entirely explicit.  Indeed, each line above may be computed using known quantities from previous lines.  However, it does not conserve energy as its Eulerian analogue does.  Indeed, following the standard steps of the energy conservation proof for implicit PIC schemes, one can observe that in fact 
\begin{equation} \label{eq:brokenEcons}
\begin{split}
    \sum_p w_p \bv_p^* \cdot \left( \bv_p^{n+1} - \bv_p^n \right) &= \Delta t \sum_p w_p \bv_p^* \cdot \bE_p^{n+1/2} \\ 
    &= \Delta t \sum_p \sum_h w_p \bv_p^* \cdot \bE_h^{n+1/2} S^h(\bx_p^{n+1/2} - \bx_h) \\
    &= \Delta t | \mathbf{h} | \sum_h \bE_h^{n+1/2} \cdot \bj_h^{*,n+1/2} \\
    &= -\Delta t | \mathbf{h} | \sum_h \nabla_h \phi_h^{n+1/2} \cdot \bj_h^{*,n+1/2} \\
    &= \Delta t | \mathbf{h} | \sum_h \phi_h^{n+1/2} \nabla_h \cdot \bj_h^{*,n+1/2} \\
    &= | \mathbf{h} | \sum_h \phi_h^{n+1/2} \left( \nabla^2_h \phi_h^{n+1} - \nabla^2_h \phi_h^{n} \right) \\
    &= -\frac{1}{2} | \mathbf{h} | \sum_h \left\{ \left\| \bE_h^{n+1} \right\|^2 - \left\| \bE_h\right\|^2 \right\},
\end{split}
\end{equation}
where we have used the integration-by-parts identities \eqref{eq:IBP_ES} that the spatial discretization is assumed to satisfy.  

While $\bv_p^*$ is an estimate of the half-step velocity, it is not exactly equal to the mean of $\bv_p^n$ and $\bv_p^{n+1}$, and thus
\begin{equation} \label{eq:fielddefs}
    \bv_p^* \cdot \left( \bv_p^{n+1} - \bv_p^n \right) \neq \frac{1}{2} \left( \left\| \bv_p^{n+1} \right\|^2 - \left\| \bv_p^n \right\|^2 \right).
\end{equation}
One can thus not interpret the very first expression in \eqref{eq:brokenEcons} as a difference of total kinetic energies, and we see that the change in potential energy is \textit{not} the negation of change in kinetic energy.  The core contribution of this paper is to propose a simple correction that remedies this.  

We modify the final line of the PIC scheme above to be the following instead:
\begin{equation} \label{eq:PICmodification}
\begin{split}
        \bv_p^\dagger &= \bv_p^n + \Delta t \left( \bE^{n+1/2}_p + \bv_p^* \times \bB(\bx_p^*, t^{n+1/2}) \right), \\
        \bv_p^{n+1} &= \mathbf{G}(\bv_p^n, \bv_p^*, \bv_p^\dagger),
\end{split}
\end{equation}
where $\mathbf{G}$ is some momentarily unspecified function.  We wish to choose $\mathbf{G}$ in such a way that
\begin{equation} \label{eq:constraint}
    \bv_p^* \cdot \left( \bv_p^\dagger - \bv_p^n \right) = \frac{1}{2} \left( \left\| \bv_p^{n+1} \right\|^2 - \left\| \bv_p^n \right\|^2 \right)
\end{equation}
holds.  Doing so will ensure energy conservation, since we will then have
\begin{equation} \label{eq:exp_econs_simple}
\begin{split}
    \frac{1}{2} \sum_p w_p \left( \| \bv_p^{n+1} \|^2 - \| \bv_p^n \|^2 \right) &= \sum_p w_p \bv_p^* \cdot \left( \bv_p^{\dagger} - \bv_p^n \right) \\
    &= \Delta t \sum_p w_p \bv_p^* \cdot \bE_p^{n+1/2} \\ 
    &= -\frac{1}{2} | \mathbf{h} | \sum_h \left\{ \left\| \bE_h^{n+1} \right\|^2 - \left\| \bE_h \right\|^2 \right\},
\end{split}
\end{equation}
using exactly the same logic applied above.  

However, \eqref{eq:constraint} is a scalar constraint on a vector quantity, so one expects there to be many functions $\mathbf{G}$ satisfying this constraint. We of course also wish to find a $\mathbf{G}$ that results in a second-order accurate scheme, and we use this freedom to ensure that.  We already have a second-order accurate estimate of the updated velocity in hand -- namely, $\bv_p^\dagger$.  It is thus sensible to define $\mathbf{G}$ such that
\begin{equation}
    \bv_p^{n+1} = \argmin_{\bv} \left\| \bv - \bv_p^\dagger \right\|^2 \quad \text{s.t.} \quad \bv_p^* \cdot \left( \bv_p^\dagger - \bv_p^n \right) = \frac{1}{2} \left( \left\| \bv \right\|^2 - \left\| \bv_p^n \right\|^2 \right).
\end{equation}
One hopes that upon solving this constrained optimization problem, we will have $\| \bv_p^{n+1} - \bv_p^\dagger \| = \mathcal{O}\left( \Delta t^3 \right)$, thus guaranteeing that we retain second-order accuracy.  

Happily, this constrained minimization problem has an analytic solution, which may be found via standard Lagrange multiplier techniques.  The result is
\begin{equation} \label{eq:opt_solution}
    \bv^{n+1}_p = \Gamma_p^n \bv^\dagger_p, \qquad \Gamma_p^n = \sqrt{ 1 + 2\frac{ \left( \bv_p^* - \frac{\bv_p^\dagger + \bv_p^n}{2} \right) \cdot (\bv_p^\dagger - \bv_p^n) }{\| \bv_p^\dagger \|^2}}.
\end{equation}
It only remains to verify that the scheme remains second-order.  To do so, it suffices to show that $\Gamma_p^n = 1 + \mathcal{O}\left( \Delta t^3 \right)$.  Using the definitions of the temporal updates in \eqref{eq:naivetranslation} and \eqref{eq:PICmodification}, one finds
\begin{equation} \label{eq:Gamma_manimputlaion}
\begin{split}
    \Gamma_p^n = \left( 1 + \Delta t^2 \frac{ \left( \bE_p^{n,*} - \bE_p^{n+1/2} \right) \cdot \left( \bE_p^{n+1/2} + \bv_p^* \times \bB(\bx_p^*, t^{n+1/2}) \right) }{\left\| \bv_p^\dagger \right\|^2} \right)^{1/2}.
\end{split}
\end{equation}
Next, note that
\begin{equation} \label{eq:EdiffAnalysis}
\begin{split}
    \bE_p^{n,*} - \bE_p^{n+1/2} &= \underbrace{\sum_h \left[ \bE_h^n - \bE_h^{n+1/2} \right] S^h(\bx_h - \bx_p^*)}_{:=\epsilon_1} \\
    &\quad + \underbrace{\sum_h \bE_h^{n+1/2} \left[ S^h(\bx_h - \bx_p^*) - S^h(\bx_h - \bx_p^{n+1/2}) \right]}_{:=\epsilon_2}.
\end{split}
\end{equation}
Using the definition of the electric field update, we immediately have
\begin{equation}
    \epsilon_1 = \frac{\Delta t}{2} \sum_h \nabla_h (\nabla_h^2)^{-1} \left( \nabla_h \cdot \bj_h^{*,n+1/2} \right) S^h(\bx_h - \bx_p^*),
\end{equation}
which is $\mathcal{O}(\Delta t)$ by the boundedness of current density and the shape function $S^h$.  

To treat $\epsilon_2$, we work in the small $\Delta t$ limit in which we may assume $\bx_p^*$ and $\bx_p^{n+1/2}$ lie in the same cell, so that $S^h$ is smooth on the line segment connecting them and may be Taylor expanded.  We Taylor expand both instances of $S^h$ about $\bx_h - (\bx_p^* + \bx_p\nph)/2$ to find, upon noting that second-order terms in the expansion cancel,  
\begin{equation} \label{eq:eps2}
\begin{split}
    \epsilon_2 &= -\sum_h \bE_h^{n+1/2} \nabla S^h(\bx_h - (\bx_p^* + \bx_p^{n+1/2})/2) \cdot \left( \bx_p^* - \bx_p^{n+1/2} \right) \\
    &\qquad + \mathcal{O} \left( \left\| \bx_p^* - \bx_p^{n+1/2} \right\|^3 \right) \\
    &= \left( \bx_p^* - \bx_p^{n+1/2} \right) \cdot \left. \nabla_\bx \left\{ \sum_h \bE_h^{n+1/2} S^h(\bx_h - \bx) \right\} \right\rvert_{\bx = \left( \bx_p^* + \bx_p^{n+1/2} \right)/2} \\
    &\qquad + \mathcal{O} \left( \left\| \bx_p^* - \bx_p^{n+1/2} \right\|^3 \right)
\end{split}
\end{equation}
By the particle update definition,
\begin{equation}
    \bx_p^* - \bx_p^{n+1/2} = \frac{\Delta t}{2} \left( \bv_p^n - \bv_p^* \right) = \order{\Dt^2}.
\end{equation}
The term in the gradient in the last expression in \eqref{eq:eps2} is precisely the continuous representation of the electric field induced by the shape function $S^h$.  This is bounded for reasonable $S^h$ and electric fields with bounded gradients, so we conclude that $\epsilon_2 = \mathcal{O}(\Delta t^2)$.  

As a result of the analysis of both $\epsilon_1$ and $\epsilon_2$, we have $\bE_p^{n,*} - \bE_p^{n+1/2} = \order{\Dt}$. Substituting this into the expression for $\Gamma_p^n$ appearing in \eqref{eq:Gamma_manimputlaion} and Taylor expanding confirms that, indeed, $\Gamma_p^n = 1 + \order{\Dt^3}$.

We note that the computation of $\Gamma_p^n$ is local to each particle, using only previous stage information, and is thus readily parallelizable and no more computationally demanding than an additional explicit stage in the particle push.  However, it must be noted that there is no guarantee that $\Gamma_p^n$ is real for every particle and every time-step.  This is a consequence of the constraint in the optimization problem, which may be rearranged to read $\| \bv \|^2 = \| \bv_p^n \|^2 + 2\bv_p^* \cdot \left( \bv_p^\dagger - \bv_p^n \right)$.  Clearly, the right-hand side of this expression must be non-negative for this constraint to be realizable by a real vector $\bv$.  We can guarantee this is so in the limit $\Delta t \rightarrow 0$, since the second term is $\mathcal{O}(\Delta t)$, but for any finite $\Delta t$ there may be particles for which $\| \bv_p^n \|^2 + 2\bv_p^* \cdot \left( \bv_p^\dagger - \bv_p^n \right) < 0$, which will lead to imaginary values of $\Gamma_p^n$.  Indeed, the expression from $\Gamma_p^n$ may be rearranged to read
\begin{equation}
    \Gamma_p^n = \frac{ \sqrt{\left\| \bv_p^n \right\|^2 +2 \bv_p^* \cdot \left( \bv_p^\dagger - \bv_p^n \right)}}{\left\| \bv_p^\dagger \right\|},
\end{equation}
which makes this fact immediately evident.  

In practice, we find that particles resulting in imaginary values of $\Gamma_p^n$ are exceedingly rare for time-steps of practical size.  We thus recover excellent energy conservation by simply replacing $\Gamma_p^n$ with one for those few particles.  However, we show in the next subsection that adding an additional stage to the field update can make this unforunate circumstance even rarer, resulting in energy errors near double precision roundoff errors in our numerical examples.  

Another important note involves work done by the magnetic field.  As stated in Section \ref{sec:Boris}, it is desirable to retain the property that $\bB$ does identically zero work.  In the absence of an electric field, the definition of $\bv_p^\dagger$ immediately implies 
\begin{equation}
    \bv_p^* \cdot \left( \bv_p^\dagger - \bv_p^n \right) = \Dt \bv_p^* \cdot \left( \bv_p^* \times \bB(\bx_p^*, t^{n+1/2} \right) = 0.
\end{equation}
Thus, by construction, the constraint \eqref{eq:constraint} assures us that $\left\| \bv_p^{n+1} \right\|^2 = \left\| \bv_p^{n} \right\|^2$ when the electric field vanishes.  Thus, the magnetic field can do no work, as in the Boris scheme.  

This is true independent of the details of the defintion of $\bv_p^*$.  The choice we've made is motivated by the fact that, when $\bE=0$, the entire particle push becomes
\begin{equation}
\begin{split}
    \bx_p^* &= \bx_p^n + \frac{\Delta t}{2} \bv_p^n, \\
    \bv_p^* &= \bv_p^n + \frac{\Delta t}{2} \bv_p^* \times \bB(\bx_p^*, t^{n+1/2}), \\
    \bv_p^{n+1} &= \bv_p^n + \Delta t \bv_p^* \times \bB(\bx_p^*, t^{n+1/2}), \\
    \bx_p^{n+1} &= \bx_p^n + \Delta t \bv_p^*.
\end{split}
\end{equation}
The middle two lines imply that $\bv_p^{n+1/2} = \bv_p^*$, and a simple rewriting of the final line reveals that this is \textit{identical} to the symmetric Boris scheme \eqref{eq:symm_Boris} in the absence of $\bE$.  Thus, this particle push treats the magnetic field in exactly the same manner as (symmetric) Boris, and only differs in its treatment of $\bE$.

\subsection{$\Gamma_p^n$ correction}
\label{sec:correction}
Observe that in \eqref{eq:Gamma_manimputlaion}, the third factor of $\Delta t$ comes from the fact that $\bE_p^{n,*}$ and $\bE_p^{n+1/2}$ approximate the electric field at points in time separated by $\Delta t/2$.  If we could arrange the scheme such that the two different fields appearing in this expression are different approximations of $\bE$ at \textit{the same time}, one might expect to get an additional factor of $\Delta t$ which results in values of $\Gamma_p^n$ that are even closer to unity and thus less likely to be imaginary.  This can in fact be achieved straightforwardly with the following modification to the scheme in the previous section: 
\begin{equation} \label{eq:good_translation}
\begin{split}
    \bx^*_p &= \bx^n_p + \frac{\Delta t}{2} \bv^n_p, \\
    \nabla_h^2 \phi_h^* &= \nabla_h^2 \phi_h^n + \frac{\Delta t}{2} \nabla_h \cdot \bj_h^{n,*}, \\
    \bv^*_p &= \bv^n_p + \frac{\Delta t}{2} \left( \bE^{*,*}_p + \bv_p^* \times \bB(\bx_p^*, t^{n+1/2})\right), \\
    \bx^{n+1}_p &= \bx^n_p + \Delta t \bv^*_p, \\
	\nabla_h^2 \phi_h^{n+1} &= \nabla_h^2 \phi_h^n + \Dt \nabla_h \cdot \bj^{*,n+1/2}_h, \\
    \bv^{\dagger}_p &= \bv^n_p + \Delta t \left(\bE^{n+1/2}_p + \bv_p^* \times \bB(\bx_p^{n+1/2}, t^{n+1/2})\right), \\
    \bv^{n+1}_p &= \bv_p^\dagger \sqrt{ 1 + 2\frac{ \left( \bv_p^* - \frac{\bv_p^\dagger + \bv_p^n}{2} \right) \cdot (\bv_p^\dagger - \bv_p^n) }{\| \bv_p^\dagger \|^2}}.
\end{split}
\end{equation}
Note that we now omit definitions that are readily understood from context.  We simply clarify that
\begin{equation} \label{eq:doublestardef}
    \bj_h^{n,*} = \frac{1}{|\mathbf{h}|} \sum_p w_p \bv_p^n S^h(\bx_h - \bx_p^*), \qquad \bE_p^{*,*} = \sum_h \bE_h^* S^h(\bx_h - \bx_p^*).
\end{equation}
With this simple modification, it is clear that \eqref{eq:Gamma_manimputlaion} becomes
\begin{equation}
    \Gamma_p^n = \left( 1 + \Delta t^2 \frac{ \left( \bE_p^{*,*} - \bE_p^{n+1/2} \right) \cdot \left( \bE_p^{n+1/2} + \bv_p^* \times \bB(\bx_p^*, t^{n+1/2}) \right)}{\left\| \bv_p^\dagger \right\|^2} \right)^{1/2}.
\end{equation}
Similar, albeit more technical, techniques to those used in the section above can be used to show that this results in $\Gamma_p^n = 1 + \mathcal{O}(\Delta t^4)$.  The derivation is sufficiently cumbersome that we place it in Appendix A.  

Thus, for a small but finite $\Delta t$, we expect it now to be even rarer that $\Gamma_p^n$ is imaginary.  This is born out in our numerical examples below.  We do note, however, that this correction does come with an additional computational cost.  An additional current deposition is required to compute $\bj_h^{n,*}$, and deposition is frequently among the rate-limiting steps in PIC algorithms.  We leave it to users to determine whether the additional energy accuracy admitted by this correction is worth this cost for their particular problem of interest.

\subsection{Straightforward electromagnetic extension}
The most natural extension of the scheme above to the electromagnetic case maintains the analogy with the motivating Eulerian scheme \eqref{eq:eulerian_scheme} by using a linearly implicit solve for the vacuum portion of Maxwell's equations: 
\begin{equation} \label{eq:EMbasic}
\begin{split}
    \bx^*_p &= \bx^n_p + \frac{\Delta t}{2} \bv^n_p, \\
    \bE_h^* &= \bE_h^n + \frac{\Delta t}{2} \left( c^2 \nabla_h \times \bB_h^n -  \bj_h^{n,*} \right), \\
    \bB_h^* &= \bB_h^n - \frac{\Delta t}{2} \nabla_h \times \bE_h^n, \\
    \bv^*_p &= \bv^n_p + \frac{\Delta t}{2} \left( \bE^{*,*}_p + \bv_p^* \times \bB_p^{*,*} \right), \\
	\bE_h^{n+1} &= \bE_h^n + \Delta t \left( c^2 \nabla_h \times \bB_h^{n+1/2} -  \bj_h^{*,n+1/2} \right), \\
    \bB_h^{n+1} &= \bB_h^n - \Delta t \nabla_h \times \bE_h^{n+1/2}, \\
    \bx^{n+1}_p &= \bx^n_p + \Delta t \bv^*_p, \\
    \bv^{\dagger}_p &= \bv^n_p + \Delta t \left( \bE^{n+1/2}_p + \bv_p^* \times \bB_p^{n+1/2} \right), \\
    \bv^{n+1}_p &= \bv_p^\dagger \sqrt{ 1 + 2\frac{ \left( \bv_p^* - \frac{\bv_p^\dagger + \bv_p^n}{2} \right) \cdot (\bv_p^\dagger - \bv_p^n) }{\| \bv_p^\dagger \|^2}}.
\end{split}
\end{equation}
Note that $\bB_p^{n+1/2}$ and $\bB_p^{*,*}$ are defined in precisely the same manner as $\bE_p^{n+1/2}$ and $\bE_p^{*,*}$ in \eqref{eq:fielddefs} and \eqref{eq:doublestardef}, respectively.

Note that the temporal discretization of Maxwell's equations is simply the Crank-Nicolson (CN) scheme.  Thus, like the semi-implicit schemes of Lapenta cited in the introduction, this scheme is implicit in the field solve while remaining explicit in all particle variables.  

By construction, we continue to have
\begin{equation} \label{eq:still_e_dot_j}
    \frac{1}{2} \sum_p w_p \left( \| \bv_p^{n+1} \|^2 - \| \bv_p^n \|^2 \right) = \Delta t | \mathbf{h} | \sum_h \bE_h^{n+1/2} \cdot \bj_h^{*,n+1/2}.
\end{equation}
The linearly implicit discretization of Maxwell's equations further implies
\begin{equation} \label{eq:EM_econs}
\begin{split}
    \Delta t \sum_h \bE_h^{n+1/2} \cdot \bj_h^{*,n+1/2} &= -\sum_h\bE_h^{n+1/2} \cdot \left( \bE_h^{n+1} - \bE_h^n - \Delta t c^2 \nabla_h \times \bB_h^{n+1/2} \right) \\
    &= -\frac{1}{2}\sum_h \left( \| \bE_h^{n+1} \|^2 - \| \bE_h^n \|^2 \right) \\
    &\qquad + \Delta t c^2 \sum_h \bE_h^{n+1/2} \cdot \nabla_h \times \bB_h^{n+1/2} \\
    &= -\frac{1}{2}\sum_h \left( \| \bE_h^{n+1} \|^2 - \| \bE_h^n \|^2 \right) \\
    &\qquad + \Delta t c^2 \sum_h \bB_h^{n+1/2} \cdot \nabla_h \times \bE_h^{n+1/2} \\
    &= -\frac{1}{2}\sum_h \left( \| \bE_h^{n+1} \|^2 - \| \bE_h^n \|^2 \right) \\
    &\qquad - c^2 \sum_h \bB_h^{n+1/2} \cdot \left( \bB_h^{n+1} - \bB_h^n \right) \\
    &= -\frac{1}{2}\sum_h \left( \| \bE_h^{n+1} \|^2 + c^2 \| \bB_h^{n+1} \|^2 - \| \bE_h^n \|^2 - c^2 \| \bB_h^n \|^2 \right).
\end{split}
\end{equation}
By combining this with \eqref{eq:still_e_dot_j}, one easily observes energy conservation.  

We note that \eqref{eq:EM_econs} relies on the identity 
\begin{equation} \label{eq:curlIdentity}
    \sum_h \left( \mathbf{G}_h \cdot \nabla_h \times \mathbf{F}_h - \mathbf{F}_h \cdot \nabla_h \times \mathbf{G}_h \right) = 0
\end{equation}
for arbitrary smooth vector fields $\mathbf{F}$ and $\mathbf{G}$ satisfying appropriate boundary conditions, which is of course true in the continuum but should be preserved by the spatial discretization for the proof above to carry through.  The Yee lattice \cite{Yee66} is well known to preserve this identity, and we show in Appendix B that this identity holds for pseudospectral discretization as well.  

Additionally, the logic in previous sections carries through with no meaningful changes to find that we still have $\Gamma_p^n = 1 + \mathcal{O}(\Delta t^4)$.  


\subsection{Charge Conservation} \label{sec:chargeconservation}
Recall from Section \ref{sec:chrgcons_bckgrnd} that by charge conservation we mean satisfaction of a discrete continuity equation.  Here, the relevant current used in updating Maxwell's equations is $\bj_h^{*,n+1/2}$, so that is the current that should also appear in the continuity equation if we are to connect the Maxwell update and continuity to realize enforcement of Gauss' law.  

In particular, we seek a definition of $\rho_h^n$ such that
\begin{equation} \label{eq:discCont}
    \frac{\rho_h^{n+1} - \rho_h^n}{\Dt} = - \nabla_h \cdot \bj_h^{*,n+1/2},
\end{equation}
as this will imply that when applying the scheme \eqref{eq:EMbasic}, we have
\begin{equation}
\begin{split}
    \nabla_h \cdot \left( \frac{\bE_h^{n+1} - \bE_h^n}{\Dt} \right) &= \nabla_h \cdot \left( c^2\nabla_h \times \bB_h\nph - \bj_h^{*,n+1/2} \right) \\
    &= \frac{\rho_h^{n+1} - \rho_h^n}{\Dt},
\end{split}
\end{equation}
where the first line is just the divergence of the discrete electric field update, and we've again assumed our spatial discretization preserves the fact that the divergence of a curl vanishes.  This in turn implies that if Gauss' law is satisfied initially, it is satisfied at all times.  

Assuming that $\rho_h^n$ has the form
\begin{equation}
    \rho_h^n = \frac{1}{| \mathbf{h} |} \sum_p w_p \mathbb{S}^h \left( \bx_h - \bx_p^n \right),
\end{equation}
a sufficient condition for \eqref{eq:discCont} to be satisfied is 
\begin{equation} \label{eq:detbal}
    \mathbb{S}^h\left( \bx_h - \bx_p^{n+1} \right) - \mathbb{S}^h\left( \bx_h - \bx_p^{n} \right) = -\left( \bx_p^{n+1} - \bx_p^n \right) \cdot \nabla_h S^h \left( \mathbf{x}_h - \bx_p^{n+1/2} \right),
\end{equation}
where we've just used the definition of $\bj_h^{*,n+1/2}$ and the fact that $\Dt \bv_p^* = \bx_p^{n+1} - \bx_p^n$.  As noted in several other works \cite{chen2011energy, chen2020semi, ricketson2023pseudospectral}, if $\mathbb{S}^h$ is quadratic on the line segment connecting $\bx_p^n$ and $\bx_p^{n+1}$, the left side can be written exactly as the analytic derivative of $\mathbb{S}^h$ evaluated at $\bx_p^{n+1/2}$ dotted with $(\bx_p^n - \bx_p^{n+1})$.  Of critical importance here is that the $\bx_p^{n+1/2}$ appearing on the right is exactly the mean of the $\bx_p^n$ and $\bx_p^{n+1}$ appearing on the left, as \textit{centered} finite differences are exact for quadratic functions.  This is what motivated our earlier choice to evaluate $\bj$ and $\bE$ at $\bx_p^{n+1/2}$ rather than, say, $\bx_p^*$.  

As a result, we have equality if and only if $\mathbb{S}^h$ is locally quadratic \textit{and} the analytic gradient of $\mathbb{S}^h$ is equal to the numerical gradient induced by the spatial discretization applied to $S^h$.
Shape function pairs $(S^h, \mathbb{S}^h)$ that satisfy these constraints have been derived in a variety of contexts.  Uniformly, they require breaking particle trajectories into segments that lie within individual cells when performing current deposition and electric field interpolation, as this is the only way to ensure that shape functions are smooth on the line segments connecting the required spatial points.  The manner in which this is done is outlined in each of the references above, as is the straightforward extension of the energy conservation proof that is necessitated by this decomposition.  That extension works in identical fashion for our scheme, but is quite lengthy, so we choose not to reprint it here.  

In one dimension with second-order finite differences, it was shown in \cite{chen2011energy} that if $S^h$ is a B-spline of degree $p$, then $\mathbb{S}^h$ is a B-spline of degree $p+1$.  Note that one requires $p \in \{ 0, 1 \}$ to guarantee that $\mathbb{S}^h$ is quadratic.  The extension to multiple dimensions with the Yee lattice is performed in \cite{chen2020semi}, where it was found that different shape functions $S^h$ must be used for each component of $\bj$, mixing degree zero and one B-splines.  Shape function pairs for pseudospectral spatial discretizations were found in \cite{ricketson2023pseudospectral} -- see particularly Appendix B in that paper.  There, the shape functions $\mathbb{S}^h$ for $\rho$ do not have compact support, so it is more computationally efficient to compute $\rho$ using the continuity equation rather than direct particle deposition.  

In each case, the fundamental starting point is \eqref{eq:detbal}.  Thus, all of these results apply trivially to the scheme presented here.  It is in this sense that we mean that our new explicit scheme permits exact charge conservation.


\section{Considerations of light-wave dispersion}
While the linearly implicit Crank-Nicolson (CN) temporal discretization of Maxwell's equations present in \eqref{eq:EMbasic} is the most straightforward one, and most directly analogous with the motivating Eulerian scheme \eqref{eq:eulerian_scheme}, it is by no means the only possibility.  Indeed, other temporal discretizations are often preferred in the literature due to their improved ability to capture the analytic dispersion relation of light waves.  In this section, we show that two of the more popular discretizations of Maxwell's equations are also compatible with our scheme, and that exact energy conservation can be recovered with either.  


\subsection{Leapfrog} \label{sec:LeapfrogMaxwell}
The combination of the Yee lattice with a leapfrog-type temporal discretization of the electromagnetic fields has been a \textit{defacto} standard in electromagnetic PIC simulation for many years.  It is often called the finite difference time domain (FDTD) method.  This may be incorporated into our development here by modifying \eqref{eq:EMbasic} as follows: 
\begin{equation} \label{eq:EMleapfrog}
\begin{split}
    \bx^*_p &= \bx^n_p + \frac{\Delta t}{2} \bv^n_p, \\
    \bB_h^{n+1/2} &= \bB_h^{n-1/2} - \Delta t \nabla_h \times \bE_h^{n}, \\
    \bE_h^* &= \bE_h^n + \frac{\Delta t}{2} \left( c^2 \nabla_h \times \bB_h^{n+1/2} -  \bj_h^{n,*} \right), \\
    \bv^*_p &= \bv^n_p + \frac{\Delta t}{2} \left( \bE^{*,*}_p + \bv_p^* \times \bB_p^{n+1/2,*} \right), \\
	\bE_h^{n+1} &= \bE_h^n + \Delta t \left( c^2 \nabla_h \times \bB_h^{n+1/2} -  \bj_h^{*,n+1/2} \right), \\
    \bx^{n+1}_p &= \bx^n_p + \Delta t \bv^*_p, \\
    \bv^{\dagger}_p &= \bv^n_p + \Delta t \left( \bE^{n+1/2}_p + \bv_p^* \times \bB_p^{n+1/2} \right), \\
    \bv^{n+1}_p &= \bv_p^\dagger \sqrt{ 1 + 2\frac{ \left( \bv_p^* - \frac{\bv_p^\dagger + \bv_p^n}{2} \right) \cdot (\bv_p^\dagger - \bv_p^n) }{\| \bv_p^\dagger \|^2}}.
\end{split}
\end{equation}
Discrete differential operators are understood here to represent the usual finite differences on the Yee lattice, so that \eqref{eq:curlIdentity} still holds.  Note that the ordering of the equations has been modified now that the updated magnetic field $\bB_h^{n+1/2}$ can be computed before computing the updated electric field, and thus before the updated velocity.  

Having made no modification to the particle push or electric field update, we still have \eqref{eq:still_e_dot_j}.  Following the line of analysis in \eqref{eq:EM_econs}, we arrive at
\begin{equation}
\begin{split}
    \Delta t \sum_h \bE_h^{n+1/2} \cdot \bj_h^{*,n+1/2} &= -\frac{1}{2} \sum_h \left( \left\| \bE_h^{n+1} \right\|^2 - \left\| \bE_h^n \right\|^2 \right) \\
    &\qquad + \Delta t c^2 \sum_h \bB_h^{n+1/2} \cdot \nabla_h \times \bE_h^{n+1/2}.
\end{split}
\end{equation}
Following exactly the analysis in \cite{chen2020semi}, we have that 
\begin{equation}
\begin{split}
    \Delta t \nabla_h \times \bE_h\nph &= \frac{1}{2} \Delta t \left( \nabla_h \times \bE_h^n + \nabla_h \times \bE_h^{n+1} \right) \\
    &= -\frac{1}{2} \left( \bB_h^{n+1/2} - \bB_h^{n-1/2} + \bB_h^{n+3/2} - \bB_h^{n+1/2} \right) \\
    &= -\frac{1}{2} \Delta t \left( \bB_h^{n+3/2} - \bB_h^{n-1/2} \right).  
\end{split}
\end{equation}

Defining electric and magnetic potential energy as 
\begin{equation}
    W_E^n = \frac{1}{2} |\mathbf{h}| \sum_h \left\| \bE_h^n \right\|^2, \qquad W_B^n = \frac{c^2}{2} |\mathbf{h}| \sum_h \bB_h^{n+1/2} \cdot \bB_h^{n-1/2},
\end{equation}
we see that 
\begin{equation}
    \Delta t | \mathbf{h} | \sum_h \bE_h^{n+1/2} \cdot \bj_h^{*,n+1/2} = - \left( W_E^{n+1} + W_B^{n+1} - W_E^n - W_B^n \right),
\end{equation}
and exact energy conservation is again recovered by combining with \eqref{eq:still_e_dot_j}.  

As discussed in \cite{chen2020semi}, the definition of magnetic potential is non-standard.  However, it is shown there that it is a second-order (in time) approximation of the more standard definition and almost always well-posed (in the sense of being non-negative).  As such, the scheme here is compatible with the FDTD discretization of the field equations in the same sense as that in \cite{chen2020semi}.


\subsection{Pseudospectral Analytic Time Domain (PSATD)}
Recall from Section \ref{sec:space} that \eqref{eq:PSATDformula} gives (the Fourier transform of) an exact solution of the ODE \eqref{eq:PSATDconstruction} for an arbitrary, fixed current density.  This analytic solution is available at any time in the interval $[t^n, t^{n+1}]$.  If we use the current density $\bj^{*,n+1/2}$ as we have in the rest of our scheme's development, then dotting the electric field update equation with $\bE_h$ and integrating in time from $t^n$ to $t^{n+1}$ gives
\begin{equation}
\begin{split}
    \bj_h^{*,n+1/2} \cdot \int_{t^n}^{t^{n+1}} \bE_h(t) \rd{t} &= - \frac{1}{2} \left( \left\| \bE_h^{n+1} \right\|^2 - \left\| \bE_h^{n} \right\|^2 \right) \\
    &\qquad + c^2 \int_{t^n}^{t^{n+1}} \bE_h \cdot \left( \nabla_h \times \bB_h \right) \rd{t}.
\end{split}
\end{equation}

For brevity, let us define 
\begin{equation}
    \left\langle \bE_h \right\rangle_n^{n+1} \equiv \frac{1}{\Dt} \int_{t^n}^{t^{n+1}} \bE_h(t)\rd{t}.
\end{equation}
Then, by summing over $h$, using the identity \eqref{eq:curlIdentity} and the magnetic field update equation dotted with $\bB_h$, we conclude that
\begin{equation} \label{eq:PSATDecons}
    \Dt \sum_h \bj_h^{*,n+1/2} \cdot \left\langle \bE_h \right\rangle_n^{n+1} = -\frac{1}{2} \sum_h \left( \left\| \bE_h^{n+1} \right\|^2 + c^2 \left\| \bB_h^{n+1} \right\|^2 - \left\| \bE_h^{n} \right\|^2 - c^2 \left\| \bB_h^{n} \right\|^2 \right).
\end{equation}

This immediately tells us how our scheme should be modified to accommodate the PSATD field solve while retaining energy conservation.  Again for brevity, denote the PSATD field update defined by \eqref{eq:PSATDformula} by $(\bE_h^{n+1}, \bB_h^{n+1}) = \text{PSATD} (\bE_h^n, \bB_h^n, \bj(\tilde{t}), \Dt )$.  We then write the PIC scheme
\begin{equation} \label{eq:EM_PSATD}
\begin{split}
    \bx^*_p &= \bx^n_p + \frac{\Delta t}{2} \bv^n_p, \\
    \left( \bE_h^*, \bB_h^* \right) &= \text{PSATD}\left( \bE_h^n, \bB_h^n, \bj_h^{n,*}, \Dt/2 \right), \\
    \bv^*_p &= \bv^n_p + \frac{\Delta t}{2} \left( \bE^{*,*}_p + \bv_p^* \times \bB_p^{*,*} \right), \\
    \bx^{n+1}_p &= \bx^n_p + \Delta t \bv^*_p, \\
	\left( \bE_h^{n+1}, \bB_h^{n+1} \right) &= \text{PSATD} \left( \bE_h^n, \bB_h^n, \bj_h^{*, n+1/2}, \Dt \right) \\
    \bv^{\dagger}_p &= \bv^n_p + \Delta t \left( \left\langle \bE_p \right\rangle_n^{n+1} + \bv_p^* \times \bB_p^{n+1/2} \right), \\
    \bv^{n+1}_p &= \bv_p^\dagger \sqrt{ 1 + 2\frac{ \left( \bv_p^* - \frac{\bv_p^\dagger + \bv_p^n}{2} \right) \cdot (\bv_p^\dagger - \bv_p^n) }{\| \bv_p^\dagger \|^2}},
\end{split}
\end{equation}
where we've defined 
\begin{equation}
    \left\langle \bE_p \right\rangle_n^{n+1} \equiv \sum_h \left\langle \bE_h \right\rangle_n^{n+1} S^h \left( \bx_h - \bx_p^{n+1/2} \right).
\end{equation}
The only noteworthy changes relative to \eqref{eq:EMbasic} are the field solve and the usage of $\left\langle \bE_p \right\rangle_n^{n+1}$ in place of $\bE_p^{n+1/2}$ to compute $\bv_p^\dagger$.  

Just as before, for this scheme we have
\begin{equation}
\begin{split}
    \frac{1}{2} \sum_p w_p \left( \left\| \bv_p^{n+1} \right\|^2 - \left\| \bv_p^{n} \right\|^2 \right) =& \sum_p w_p \bv_p^* \cdot \left( \bv_p^\dagger - \bv_p^n \right) \\
    &= \Delta t | \mathbf{h} | \sum_h \left\langle \bE_h \right\rangle_n^{n+1} \cdot \bj^{*,n+1/2}.
\end{split}
\end{equation}
This, in concert with \eqref{eq:PSATDecons}, gives exact energy conservation.  Trivial modifications of earlier arguments show that we still have $\Gamma_p^n = 1 + \order{\Dt^4}$, and the scheme thus remains second-order accurate in time.  

Note that an explicit form for the Fourier transform of $\left\langle \bE_h \right\rangle_n^{n+1}$ can be found by letting $\Dt \rightarrow \tau$ in \eqref{eq:PSATDformula} and integrating from $\tau=0$ to $\tau=\Dt$.  One finds
\begin{equation}
\begin{split}
    \mathcal{F} \left[ \left\langle \bE_h \right\rangle_n^{n+1} \right] &= \frac{S}{kc \Dt} \mathcal{F} \left[ \bE^n \right] + i \frac{1-C}{k\Dt} \widehat{\mathbf{k}} \times \mathcal{F} \left[ \bB^n \right] - \frac{1 - C}{k^2 c^2 \Dt} \mathcal{F} \left[ \bj^{*, n+1/2} \right] \\
    &\qquad + \left( 1 - \frac{S}{kc\Dt} \right) \widehat{\mathbf{k}} \left( \widehat{\mathbf{k}} \cdot \mathcal{F} \left[ \bE^n \right] \right) \\
    &\qquad + \widehat{\mathbf{k}} \left( \widehat{\mathbf{k}} \cdot \mathcal{F} \left[ \bj^{*,n+1/2} \right] \right) \left( \frac{1-C}{k^2 c^2 \Dt } - \frac{\Dt}{2} \right).
\end{split}
\end{equation}
Upon computing the inverse Fourier transform, we have an explicitly evaluable expression for $\left\langle \bE_h \right\rangle_n^{n+1}$, and we see that the scheme is again fully explicit.  

Note that a directly analogous line of reasoning can be used to adapt the semi-implicit scheme of \cite{chen2020semi} to the use of PSATD while keeping exact energy conservation.  One could also retain exact charge conservation by using the shape function pair introduced in \cite{ricketson2023pseudospectral}.  The testing of such a scheme is an interesting avenue for future research.  

As a final note, we acknowledge that there is some flexibility in the evaluation of $\bB$.  It may be natural to use a time-average of the analytic expression, as is done for $\bE$.  However, this is not necessary for energy conservation, and both that choice and the use of $\bB_p^{n+1/2}$ result in the same order of accuracy.  An investigation of the different properties of each choice is again an interesting topic for further research.  


\section{Numerical Results}
\label{sec:numerics}
We have implemented several of the new schemes described above in a 2D2V test code.  This code uses a pseudospectral discretization of all spatial derivative operators, accelerated as usual by the fast Fourier transform.  As such, we do not directly test the scheme in Section \ref{sec:LeapfrogMaxwell}, as this leapfrog scheme for Maxwell's equations is most logically combined with finite differences on the Yee lattice.  Additionally, our code does not enforce strict charge conservation as described in Section \ref{sec:chargeconservation}.  We rely on the theoretical results in that section as evidence that charge conservation can be enforced.  

When making comparisons, it will be useful to have names for each scheme.  To that end, we denote the initial electrostatic, energy-conserving scheme derived in Section \ref{sec:initialscheme} as ESEC1 (ElectroStatic Energy Conserving 1).  Its improvement in Section \ref{sec:correction}, written explicitly in equation \eqref{eq:good_translation}, is denoted as ESEC2.  Recall that the distinction between these two schemes is that the $\Gamma_p^n$ factor that enforces energy conservation is closer to unity for ESEC2, and thus less likely to result in imaginary values.  

For the electromagnetic case, we test the scheme with linearly implicit field solve defined in equation \eqref{eq:EMbasic} and denote it as EMEC(LI).  We also test the scheme with the PSATD field solve defined in \eqref{eq:EM_PSATD} and denote it as EMEC(PSATD).  

We recall from \cite{ricketson2023pseudospectral} that the second identity in \eqref{eq:IBP_ES} is satisfied by pseudospectral discretizations \textit{except} for the Nyquist mode.  In that work, binomial filtering is employed to eliminate the Nyquist mode and thus recover exact energy conservation.  We do the same here, noting that binomial filtering is commonly used in PIC schemes to mitigate particle sampling noise, and that the manner in which it should be used to preserve energy conservation is detailed in many of the references on implicit PIC listed above.  

In the electrostatic tests, we compare the new schemes to a standard explicit PIC scheme, which uses the Vlasov-Poisson formulation with the classical leapfrog/Boris particle push.  For consistency, we also apply binomial filtering to the computed charge densitites for the standard scheme.  

Our numerical studies are conducted in a two-dimensional periodic box with identical side lengths $L$.  We denote the number of cells in each direction as $N_x$ and $N_y$, and the number of particles per cell as $N_c$.  All shape functions used are the first-order $B$-splines, sometimes called ``hat" or ``tent" functions.  


\subsection{Linear Landau Damping}
We repeat the linear Landau damping test of \cite{ricketson2016sparse}, albeit in two dimensions rather than three.  For completeness, we reproduce the setup of that test here.  We sample particles from the initial distribution function
\begin{equation}
    f_0(x,y,v_x,v_y) = \frac{1}{2\pi} e^{-(v_x^2 + v_y^2)/2} \left( 1 + \alpha_x \cos \left( \frac{2 \pi x}{L} \right) \right) \left( 1 + \alpha_y \cos \left( \frac{2 \pi y}{L} \right) \right).
\end{equation}
For this linear damping test, we set $\alpha_x = \alpha_y = 0.05$ and use box size $L = 22$.  We set $N_x = N_y = 32$, $\Delta t = 1/10$, and use $N_c=500$ particles per cell with random sampling from $f_0$ for initialization.  

In our dimensionless variables, the linear Landau damping rate $\gamma$ for a perturbation with wave-number $k$ is given approximately by \cite{krall1973principles}
\begin{equation}
    \gamma \approx \sqrt{\frac{\pi}{2}} \frac{\omega^2}{2 k^3} e^{-\omega^2/ 2 k^2}, 
\end{equation}
with the oscillation frequency $\omega$ given in terms of $k$ by the Bohm-Gross dispersion relation
\begin{equation}
    \omega^2 = 1 + 3 k^2.
\end{equation}
In Figure 1, we compare the electrostatic potential energy evolution in ESEC1, ESEC2, and standard PIC to this theoretical damping rate.  
\begin{figure}[h]
    \centering
    \includegraphics[width=0.7\linewidth]{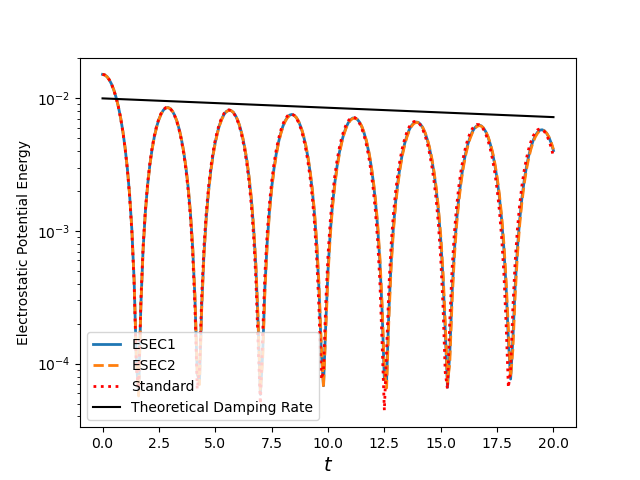}
    \caption{Evolution of electrostatic potential energy for linear Landau damping test case.  Excellent agreement is observed between all three schemes as well as the theoretical damping rate.}
    \label{fig:landau_verification}
\end{figure}
We observe excellent agreement of the new energy-conserving schemes with both standard PIC and the linear theory.  

In Figure 2, we display the fractional change in total energy over time for the three schemes tested.  Recall that total energy is defined by
\begin{equation}
    \text{TE}^n = \frac{1}{2} \sum_p w_p \left\| \bv_p^n \right\|^2 + \frac{1}{2} \lvert \mathbf{h} \rvert \sum_h \left\| \bE_h^n \right\|^2
\end{equation}
in the electrostatic case, with an extra factor of $c^2 \left\| \bB_h^n \right\|^2$ appearing in the second term for electromagnetic simulations.  We thus define fractional energy change by
\begin{equation}
    \delta^n = \frac{\text{TE}^n - \min_k \text{TE}^k}{\text{TE}^0}.
\end{equation}
\begin{figure}[h]
    \centering
    \includegraphics[width=0.7\linewidth]{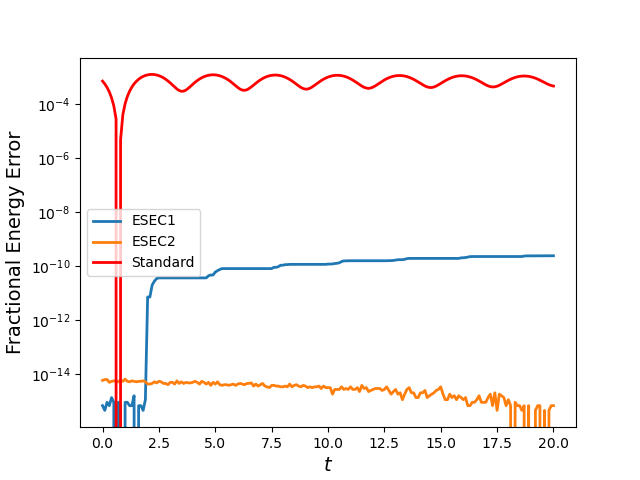}
    \caption{Fractional change in total energy over time for the three tested schemes.}
    \label{fig:econs_linearlandau}
\end{figure}
As predicted by our theory, the new energy conserving schemes feature dramatically reduced energy errors compared to standard PIC.  Indeed, ESEC2 conserves energy to machine precision.  ESEC1 only fails to reach machine precision because of a small number of particles at some time-steps for which $\Gamma_p^n$ is imaginary.  Recall that to avoid imaginary velocities, we artificially set $\Gamma_p^n = 1$ which this situation arises, thus making small energy errors.  

In Figure 3, we justify this interpretation by plotting the number of particles for which $\Gamma_p^n$ is imaginary at each time-step, both for ESEC1 and ESEC2.  
\begin{figure}[h]
    \centering
    \includegraphics[width=0.7\linewidth]{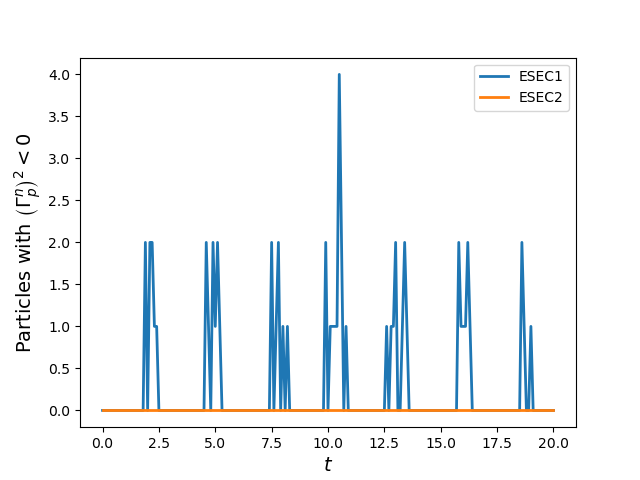}
    \caption{Number of problematic particles with imaginary correction factors $\Gamma_p^n$ at each time-step for linear Landau damping with $\Delta t = 0.1$.}
    \label{fig:enter-label}
\end{figure}
Again as predicted by our theory, such particles are rare, but even rarer for ESEC2 than ESEC1.  At this time-step size, we do not observe a single problematic particle for ESEC2.  For ESEC1, we observe at most four problematic particles at any given time-step.  Note that the total number of particles in this simulation is $N_p = 512000$, so these particles are still quite rare, leading to excellent -- but not quite exact -- energy conservation for ESEC1.  

\begin{figure}
    \centering
    \includegraphics[width=0.7\linewidth]{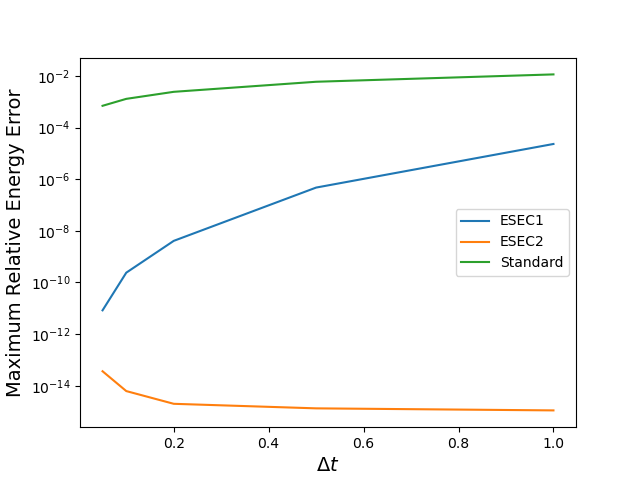}
    \caption{Maximum fractional energy errors for the linear Landau damping test as a function of time-step size.}
    \label{fig:maxerrs_vs_dt}
\end{figure}

Going further, we can study the maximum energy error for each scheme as a function of time-step size.  These results are shown in Figure \ref{fig:maxerrs_vs_dt}.
As the theory predicts, energy errors increase with time-step for ESEC1 because more particles feature imaginary values of $\Gamma_p^n$.  For this test case, no issues arise for ESEC2 even with $\Delta t = 1$ and the energy conservation is essentially exact even when the problem's dynamics are poorly resolved.  

\subsection{Two-Stream Instability}
We initialize the distribution function with two counter-streaming beams in the $x$-direction and a density perturbation in $x$:
\begin{equation}
\begin{split}
    f_0(x,y,v_x,v_y) &= \frac{1}{4\pi} e^{-v_y^2/2}\left( e^{-(v_x - v_b)^2/2} + e^{-(v_x + v_b)^2/2} \right)  \\
    &\qquad \times\left( 1 + \alpha_x \cos \left( \frac{2 \pi x}{L} \right) \right).
\end{split}
\end{equation}
We use the parameters $v_b = 3.5$, $\alpha_x = 0.01$, $L = 32$.  For discretization parameters, we choose $\Delta t = 0.1$, $N_x = N_y = 32$, and $N_c=500$ particles per cell (again giving total particle number $N_p = 512,000$).  

\begin{figure}[h]
    \centering
    \includegraphics[width=0.7\linewidth]{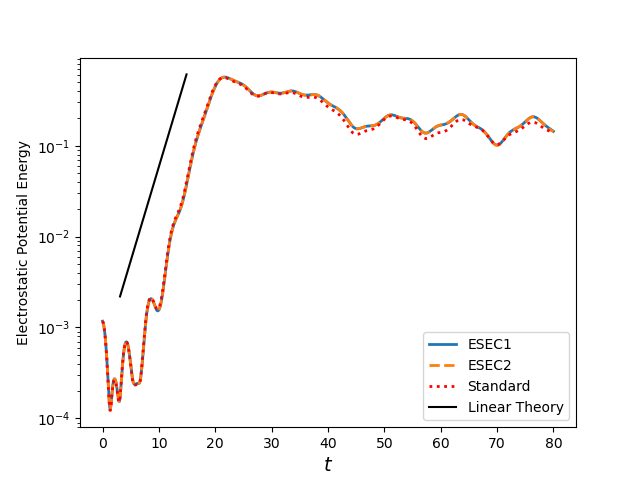}
    \caption{Electrostatic potential energy for the two-stream instability test case, showing good agreement between all schemes and the theoretical linear growth rate.}
    \label{fig:twostreamPE}
\end{figure}

The linear dispersion relation for the two-stream instability is given in our dimensionless variables by \cite{chen2011energy,stix1992waves}
\begin{equation}
    \frac{1}{(\omega + k v_b)^2} + \frac{1}{(\omega - k v_b)^2} = 1.
\end{equation}
Solving this quartic equation for $\omega$ gives two complex-conjugate solution, whose imaginary part gives the linear growth rate of the instability.  This theoretical prediction is compared against the three schemes considered here in Figure \ref{fig:twostreamPE}.
As before, we see excellent agreement.  

\begin{figure}[h]
    \centering
    \includegraphics[width=0.495\linewidth]{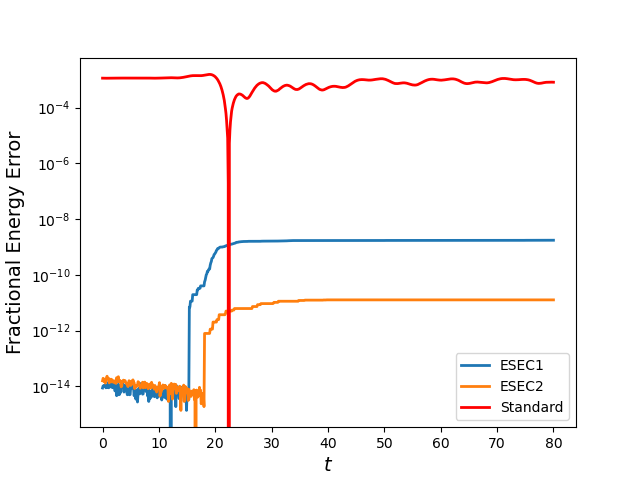}
    \includegraphics[width=0.495\linewidth]{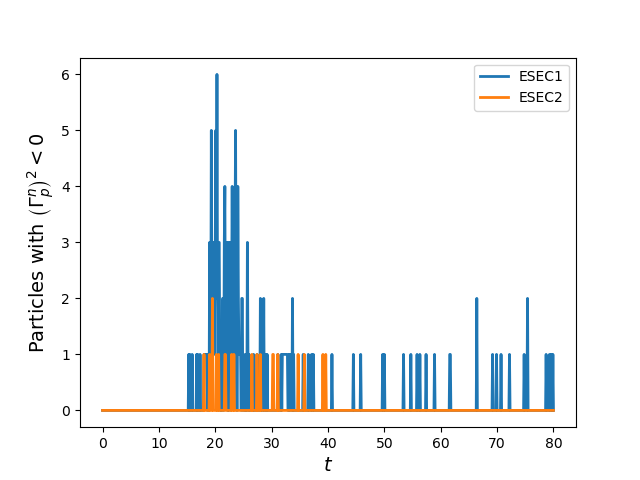}
    \caption{\underline{Left}: Fractional energy errors for the two-stream instability case as function of time for the three schemes tested.  \underline{Right}: Number of problematic particles at each time-step for ESEC1 and ESEC2.}
    \label{fig:twostreamEnergyConservation}
\end{figure}

Just as in the previous case, we plot relative energy errors in Figure \ref{fig:twostreamEnergyConservation} (left), in addition to the number of problematic particles at each time-step (right).  ESEC2 again has improved energy accuracy due to the reduced number of problematic particles.  However, in this case ESEC2 does feature a non-zero number of such particles, resulting in energy errors larger than machine precision.  
Nevertheless, energy accuracy remains improved by 6 (for ESEC1) to 8 (for ESEC2) orders of magnitude relative to standard PIC.  

We note that problematic particles are most common near the time the instability saturates.  This is likely due to fine-scale structures that develop in space near that time, which lead to large electric fields and thus rapid changes in particle velocity.  


\subsection{Weibel Instability}
We repeat the 1D2V Weibel instability test case of \cite{cheng2014energy}.  In particular, using our nondimensional formulation with $c=1$ corresponds to the scaling used there, and we initialize the distribution function and fields according to 
\begin{equation}
\begin{split}
    f(y, v_x, v_y, t=0) &= \frac{1}{\pi \beta} e^{-v_y^2/\beta} \left[ \delta e^{-(v_x - v_{0,1})^2/\beta} + (1 - \delta) e^{-(v_x + v_{0,2})^2/\beta} \right], \\
    E_x(y, t=0) &= E_y(y, t=0) = 0, \\
    B_z(y, t=0) &= b \sin (k_0 y ).
\end{split}
\end{equation}
We use the parameters from Run 1 in \cite{cheng2014energy}.  Namely, 
\begin{equation}
\begin{split}
    &\beta = 0.01, \qquad \delta = 0.5, \qquad v_{0,1} = v_{0,2} = 0.3, \\
    &b = 0.001, \qquad k_0 = 0.2.
\end{split}
\end{equation}

Because the initial perturbation amplitude (controlled by $b$) is so small, and the Weibel instability saturates at relatively low amplitude, we employ a so-called ``quiet start" for this test to observe the linear growth of the instability over several orders of magnitude.  In particular, the initial particle positions are specified deterministically on a uniform grid.  Particle velocities are still randomly sampled from the specified bi-Maxwellian distribution.    

Having thoroughly explored the differences between ESEC1 and ESEC2 in previous examples, we use this test to focus on demonstrating that both the CN 
and PSATD field solvers can produce energy conservation with our scheme.  Because the spatial discretization is necessarily spectral in our code, the Yee lattice is not tested here.  Moreover, for both PSATD and CN we always use the correction described in Section \ref{sec:correction} that has been shown in previous tests to improve energy accuracy.  
\begin{figure}[h]
    \centering
    \includegraphics[width=.495\linewidth]{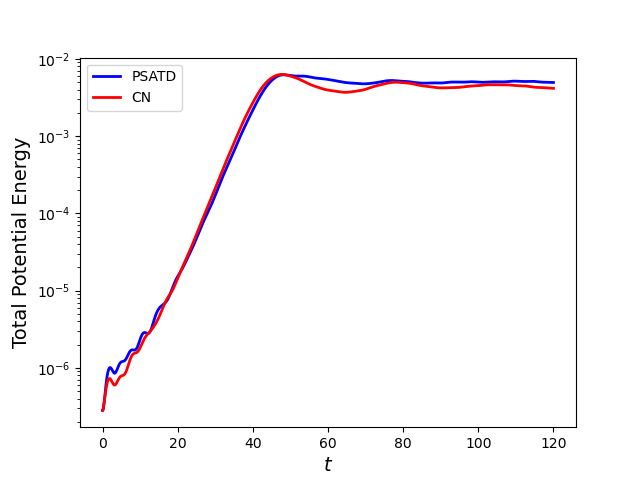}
    \includegraphics[width=.495\linewidth]{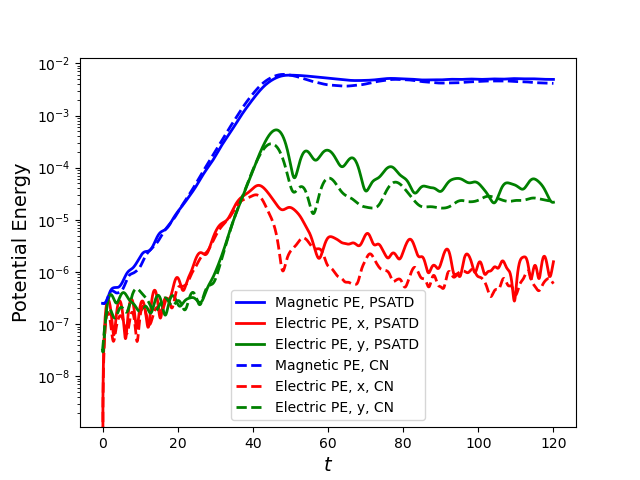}
    \caption{\underline{Left}: Total potential energy as a function of time for both the pseudo-spectral analytic time-domain (PSATD) and Crank-Nicolson (CN) temporal discretizations of the field solve.  \underline{Right}: Breakdown of various sources of potential energy, again as a function of time.}
    \label{fig:WeibelPE}
\end{figure}

We use $N_y = 32$, with $N_c=3200$ particles per cell (giving total particle number of $N_p = 1.024 \times 10^5$), domain size $L_y = 2 \pi / k_0$, and time-step $\Delta t = 0.1$.  In Figure \ref{fig:WeibelPE}, we plot potential energy -- both total and a breakdown of its various sources -- as a function of time.  We observe linear growth of the instability, as well as good agreement between the two schemes and the results in \cite{cheng2014energy}.  

\begin{figure}[h]
    \centering
    \includegraphics[width=.495\linewidth]{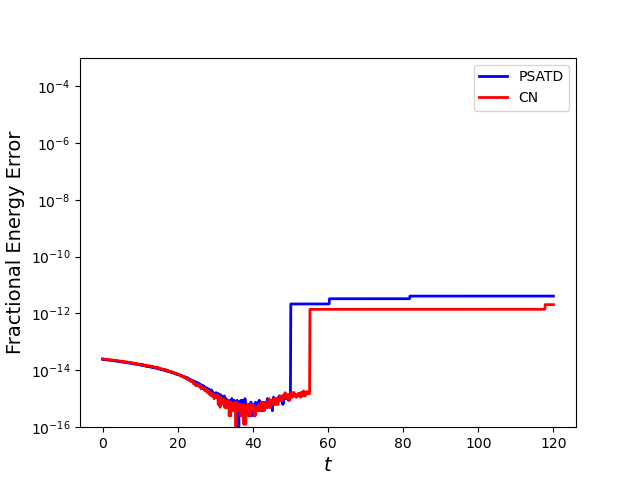}
    \includegraphics[width=.495\linewidth]{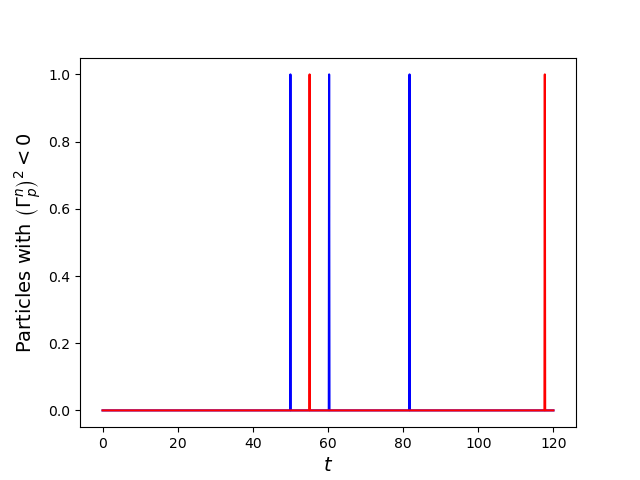}
    \caption{\underline{Left}: Fractional energy change as a function of time for both PSATD and CN schemes applied to the Weibel instability.  \underline{Right}: Problematic particles with imaginary $\Gamma_p^n$ at each time-step in the Weibel instability problem.  As with other examples, such particles are extremely rare.  }
    \label{fig:WeibelEnergy}
\end{figure}

We confirm energy conservation in Figure \ref{fig:WeibelEnergy}, which again shows roughly 11 digits of energy accuracy and that particles with imaginary $\Gamma_p^n$ are extremely rare.  

\section{Conclusions}
\label{sec:conc}
In this work, we have presented an explicit, energy-conserving particle-in-cell scheme, both for electostatic and electromagnetic cases.  The new method combines the computational efficiency and ease of implementation of existing explicit schemes with the energy conservation enjoyed by recent implicit PIC schemes.  

The local nature of the optimization procedure used to enforce energy conservation admits error analysis and is expected to facilitate scalability of the algorithm.  Additionally, we have shown that the scheme is compatible with exact local charge conservation as well as two widely used spatial discretizations known to feature desirable properties; namely, the FDTD and PSATD schemes.

There are numerous interesting directions for future work.  These include the implementation and testing of exact charge conservation schemes analyzed in theory in Section \ref{sec:chargeconservation}; extension to relativistic plasmas; studies of the scheme's behavior with respect to the finite grid instability; implementation/testing of the FDTD scheme; and many others.  

\appendix


\section*{Acknowledgements}
The authors wish to acknowledge valuable private communication with Luis Chac\'{o}n and Justin Angus.  This work was performed under the auspices of the U.S. Department of Energy by LLNL under contract DE-AC52-07NA27344. Both authors were supported by the DOE Office of Applied Scientific Computing Research (ASCR) Mathematical Multifaceted Integrated
Capabilities Center (MMICCs) Program under grant DE-SC0023164.  Additionally, the work of J. Hu was partially supported by AFOSR grant FA9550-21-1-0358.

\section*{Appendix A}
\noindent Recall from the main text that with the correction introduced in \ref{sec:correction}, we have
\begin{equation} \label{eq:GammaAppendix}
    \Gamma_p^n = \left( 1 + \Delta t^2 \frac{ \left( \bE_p^{*,*} - \bE_p^{n+1/2} \right) \cdot \left( \bE_p^{n+1/2} + \bv_p^* \times \bB(\bx_p^*, t^{n+1/2}) \right)}{\left\| \bv_p^\dagger \right\|^2} \right)^{1/2}.
\end{equation}

We analyze the difference of electric fields appearing in this expression as follows:
\begin{equation}
\begin{split}
    \bE_p^{*,*} - \bE_p^{n+1/2} &= \underbrace{\sum_h \left[ \bE^*_h - \bE_h^{n+1/2} \right] S^h(\bx_h - \bx_p^*)}_{:=\tilde{\epsilon}_1} \\
    &\quad + \underbrace{\sum_h \bE_h^{n+1/2} \left[ S^h(\bx_h - \bx_p^*) - S^h(\bx_h - \bx_p^{n+1/2}) \right]}_{:=\epsilon_2}.
\end{split}
\end{equation}
Note that $\epsilon_2$ is identical to the term of the same name analyzed in \eqref{eq:eps2}, where we concluded it was $\order{\Dt^2}$. We thus proceed with analyzing $\tilde{\epsilon}_1$.

Using the update equations for the electric field, we find
\begin{equation} \label{eq:eps1tilde}
\begin{split}
    \tilde{\epsilon}_1 &= \frac{\Delta t}{2} \sum_h \nabla_h (\nabla_h^2)^{-1} \left[ \nabla_h \cdot (\bj_h^{n,*} - \bj_h^{*,n+1/2}) \right] S^h(\bx_h - \bx_p^*).
\end{split}
\end{equation}
Next, it behooves us to consider the difference of currents
\begin{equation} \label{eq:currentdiff}
\begin{split}
     \bj_h^{n,*} - \bj_h^{*,n+1/2}  &= \frac{1}{|\mathbf{h}|}\sum_p w_p \left\{ \bv_p^n S^h(\bx_h - \bx_p^*) - \bv_p^* S^h(\bx_h - \bx_p\nph) \right\} \\
    &= \frac{1}{|\mathbf{h}|}\sum_p w_p \left( \bv_p^n - \bv_p^* \right) S^h(\bx_h - \bx_p^*) \\
    &\qquad + \frac{1}{|\mathbf{h}|} \sum_p w_p \bv_p^* \left[ S^h(\bx_h - \bx_p^*) - S^h(\bx_h - \bx_p\nph) \right] \\
    &= -\frac{\Delta t}{2} \frac{1}{|\mathbf{h}|}\sum_p \left( \bE_p^{*,*} + \bv_p^* \times \bB(\bx_p^*, t^{n+1/2}) \right) w_p S^h(\bx_h - \bx_p^*) \\
    &\qquad - \frac{1}{|\mathbf{h}|}\sum_p w_p \bv_p^* \nabla S^h \left( \bx_h - (\bx_p^* + \bx_p\nph)/2 \right) \cdot (\bx_p^* - \bx_p\nph) \\
    &\qquad + \order{\left\| \bx_p^* - \bx_p\nph \right\|^3},
\end{split}
\end{equation}
where we've used exactly the same Taylor expansion strategy as in \eqref{eq:eps2}. 
We break the first term in the last equality into two pieces, one for each term in the Lorentz force.  For the term featuring the electric field, we have
\begin{equation} \label{eq:term1}
    \frac{\Delta t}{2} \left\| \frac{1}{|\mathbf{h}|} \sum_p \bE_p^{*, *} w_p S^h(\bx_h - \bx_p^*) \right\| \leq \frac{\Delta t}{2} \left\{ \max_p \left\| \bE_p^{*,*} \right\| \right\} \left\lvert \frac{1}{|\mathbf{h}|} \sum_p w_p S^h(\bx_h - \bx_p^*) \right\rvert.
\end{equation}
Note that the last term on the right is a second-order accurate approximation of the charge density at point $\bx_h$.  Assuming both charge density and electric fields are bounded, we find that this factor is $\order{\Dt}$.  

For the term featuring the magnetic field, we have
\begin{equation}
\begin{split}
    &\frac{\Delta t}{2} \left\| \frac{1}{|\mathbf{h}|} \sum_p \bv_p^* \times \bB(\bx_p^*, t^{n+1/2}) w_p S^h(\bx_h - \bx_p^*) \right\| \\
    &\qquad \leq \frac{\Delta t}{2} \left\{ \max_p \left\| \bB(\bx_p^{*}, t^{n+1/2}) \right\| \right\} \left\lvert \frac{1}{|\mathbf{h}|} \sum_p w_p \bv_p^* S^h(\bx_h - \bx_p^*) \right\rvert.
\end{split}
\end{equation}
Much like the electric field portion, the last term on the right is a second-order approximation of the current density at point $\bx_h$.  Assuming current density and the magnetic field are both bounded again reveals that this factor is $\order{\Dt}$.  

Consider next the other term in the last equality of \eqref{eq:currentdiff}.  Using the position update definitions and introducing the name $\bx_p^+ = (\bx_p^* + \bx_p\nph)/2$ for brevity, we have 
\begin{equation} \label{eq:term2}
\begin{split}
    & \left\| \frac{1}{|\mathbf{h}|}\sum_p w_p \bv_p^* \nabla S^h \left( \bx_h - (\bx_p^* + \bx_p\nph)/2 \right) \cdot (\bx_p^* - \bx_p\nph) \right\|\\ 
    &\qquad \quad =  \left\| \frac{\Delta t^2}{4} \frac{1}{|\mathbf{h}|} \sum_p w_p \bv_p^* \left( \bE_p^{*,*} + \bv_p^* \times \bB(\bx_p^*, t\nph) \right) \nabla S^h(\bx_h - \bx_p^+) \right\| \\
\end{split}
\end{equation}
As before, we break the last line into pieces corresponding to the electric and magnetic terms in the Lorentz force.  For the electric field term, we have
\begin{equation}
\begin{split}
    &\left\| \frac{\Delta t^2}{4} \frac{1}{|\mathbf{h}|} \sum_p w_p \bv_p^* \bE_p^{*,*}  \nabla S^h(\bx_h - \bx_p^+) \right\| \\
    &\qquad \leq \frac{\Delta t^2}{4} \left\{ \max_p \left\| \bE_p^{*,*} \right\| \right\} \left\| \left. \nabla_{\bx} \left( \frac{1}{|\mathbf{h}|} \sum_p w_p \bv_p^* S^h(\bx - \bx_p^+) \right) \right\rvert_{\bx = \bx_h }\right\|.
\end{split}
\end{equation}
The last term on the right is the derivative with respect to the location of the grid point $\bx_h$ of a current density at that point.  Again assuming current density has bounded derivatives, we find that this term is indeed $\order{\Dt^2}$.  Directly analogous treatment of the magnetic field term leads to $\order{\Dt^2}$ under the assumption of bounded derivatives of the pressure tensor.  

The combination of \eqref{eq:term1} and \eqref{eq:term2} allows us to conclude that $\bj_h^{n,*} - \bj_h^{*,n+1/2} = \order{\Dt}$.  Substituting this into \eqref{eq:eps1tilde} and again assuming boundedness of spatial derivatives, we find that $\tilde{\epsilon}_1 = \order{\Dt^2}$.  This in turn implies that $\bE_p^{*,*} - \bE_p^{n+1/2} = \order{\Dt^2}$.  Substituting this into \eqref{eq:GammaAppendix} and Taylor expanding confirms that $\Gamma_p^n = 1 + \order{\Dt^4}$, as desired.

\section*{Appendix B}
As leveraged in \cite{ricketson2023pseudospectral}, the pseudospectral differentiation operator may be conceived as a matrix multiplication by a real matrix $D$ which is antisymmetric -- i.e.\ $D = -D^T$ \cite{trefethen2000spectral} .  In one dimension, the discrete analogue of integration by parts follows trivially:
\begin{equation}
    \sum_i F_i \sum_j D_{ij} G_j = -\sum_j G_j \sum_i D_{ji} F_i.
\end{equation}
Applying this component-wise to the terms in the dot products in the desired identity in three dimensions immediately gives the desired result.  For instance, for the $x$-component we have 
\begin{equation}
\begin{split}
    \sum_{ijk} F^x_{ijk} (\nabla_h \times \mathbf{G}_{ijk})^x &= \sum_{ijk} F^x_{ijk} \left( \partial_y^h G^z_{ijk} - \partial_z^h G^y_{ijk} \right) \\
    &= \sum_{ijk} F^x_{ijk} \left( \sum_l D_{jl} G^z_{ilk} - \sum_l D_{kl} G^y_{ijl} \right) \\
    &= -\sum_{ilk} G^z_{ilk} \sum_j D_{lj} F^x_{ijk} + \sum_{ijl} G^y_{ijl} \sum_k D_{lk} F_{ijk}^x \\
    &= \sum_{ijk} \left( G^y_{ijk} \partial_z^h F_{ijk}^x - G^z_{ijk} \partial_y^h F^x_{ijk} \right).
\end{split}
\end{equation}
Applying identical logic to the $y$ and $z$ components gives the desired result.  

 \bibliographystyle{elsarticle-num} 
 \bibliography{biblio}





\end{document}